\RequirePackage{fix-cm}
\documentclass[11pt, a4paper]{article}  
\providecommand{\keywords}[1]
{
  \small	
  \textbf{\textit{Keywords---}} #1
}
\newtheorem{theorem}{Theorem}
\usepackage{graphicx}
\usepackage{verbatim}
\usepackage{amssymb}
\usepackage{amsmath,amsfonts,mathtools}
\usepackage{lineno}
\usepackage[format=plain, indention=1cm]{caption}
\usepackage{subcaption}
\usepackage[breaklinks=true]{hyperref}
\usepackage{breakcites}
\usepackage{natbib}
\usepackage{lineno}
\usepackage{titles}
\usepackage{xcolor}
\usepackage{fancyhdr}
\usepackage{hyperref}

\begin{document}

\title{Multiple attractors and long transients in spatially structured populations with an Allee effect}


\author{Irina Vortkamp \and
        Sebastian J. Schreiber
        \and
        Alan Hastings
        \and
        Frank M. Hilker
}

\maketitle

\begin{abstract}
We present a discrete-time model of a spatially structured population and explore the effects of coupling when the local dynamics contain a strong Allee effect and overcompensation. While an isolated population can exhibit only bistability and essential extinction, a spatially structured population can exhibit numerous coexisting attractors. We identify mechanisms and parameter ranges that can protect the spatially structured population from essential extinction, whereas it is inevitable in the local system. In the case of weak coupling, a state where one subpopulation density lies above and the other one below the Allee threshold can prevent essential extinction. Strong coupling, on the other hand, enables both populations to persist above the Allee threshold when dynamics are (approximately) out-of-phase. In both cases, attractors have fractal basin boundaries. Outside of these parameter ranges, dispersal was not found to prevent essential extinction. We also demonstrate how spatial structure can lead to long transients of persistence before the population goes extinct.

\end{abstract}

\section{Introduction}
One of the simplest systems with the potential to exhibit a regime shift is a population with a strong Allee effect \citep{johnson2018resilience}. Population densities above a certain threshold, called Allee threshold, persist whereas populations that fall under the Allee threshold go extinct \citep{courchamp2008allee}. There is abundant evidence that Allee effects play an important role in diverse biological systems \citep{dennis1989allee, courchamp1999inverse, stephens1999allee, stephens1999consequences, courchamp2008allee}. Mechanisms that induce an Allee effect, like mate finding problems or defence against predators in small populations, are well understood \citep{courchamp2008allee}.\\
Introducing spatial structure into population models can change their dynamical behaviour.
This is of particular relevance when the local dynamics include a strong Allee effect. However, Allee effects were considered mostly in models for spatially structured populations in continuous time \citep{gruntfest1997fragmented, amarasekare1998allee, gyllenberg1999allee, kang2011expansion, wang2016population, johnson2018resilience}. One important result from these models is the rescue effect, where a subpopulation that falls under the Allee threshold is rescued from extinction by migration from another location \citep{brown1977turnover}. Moreover, \cite{amarasekare1998allee} suggests that local populations that are linked by dispersal are more abundant and less susceptible to extinction than isolated populations. Little attention has been devoted to the case in discrete time where local dynamics can be chaotic. In that case, the correlation between abundance and extinction risk 
is less obvious. There have been several studies to understand mechanisms and consequences of coupling patches in discrete time \citep{gyllenberg1993does, hastings1993complex, lloyd1995coupled, gyllenberg1996bifurcation, kendall1998spatial, earn2000coherence, yakubu2002interplay, yakubu2008asynchronous, faure2014quasi,franco2015connect}. A controversial question is whether chaotic behaviour of the population increases the probability of extinction \citep{thomas1980chaos, berryman1989ecological, lloyd1995coupled} or promotes spatially structured populations \citep{allen1993chaos} and population persistence \citep{huisman1999biodiversity}, which demands further research on coupled patches of chaotic dynamics.\\
\cite{neubert1997simple} and \cite{schreiber2003allee} study single species models with overcompensating density dependence and Allee effect. Overcompensation occurs as a lagged effect of density-dependent feedback. As a result, populations can alternate from high to low numbers even in the absence of stochasticity \citep{ranta2005ecology}. This can lead to essential extinction, a phenomenon that does not occur in corresponding continuous-time models. A major characteristic is that large population densities fall below the Allee threshold when the overcompensating response is too strong. Thus, ‘‘almost every’’ initial density leads to extinction when per capita growth is sufficiently high. In that case, \cite{schreiber2003allee} proved that long transient behaviour can occur before the population finally goes extinct. However, an interesting question that has not been studied yet is how the dynamics change when we include spatial structure. In this paper we examine the interplay between essential extinction due to local chaotic dynamics with Allee effect and the between-patch effects due to coupling.\\
We distinguish two drivers of multistability. Firstly, different states can be caused by the Allee effect \citep{dennis1989allee,gruntfest1997fragmented,amarasekare1998allee,courchamp1999inverse,gyllenberg1999allee,schreiber2003allee}. These also exist in isolated patches unless there is essential extinction. Secondly, multistability can be caused by coupling maps with overcompensation \citep{allen1993chaos,gyllenberg1993does,hastings1993complex,lloyd1995coupled,kendall1998spatial,yakubu2002interplay,wysham2008sudden,yakubu2008asynchronous}. The former occur also in continuous-time models with Allee effect, while the latter occur in discrete-time overcompensatory models without Allee effect. By including discrete-time overcompensation and Allee effects, we help to unify these separate areas of prior work.\\
The remainder of the paper is organized as follows. In Section \ref{Model}, we present an overview of the model and our main assumptions. With the aid of numerical simulations, we describe the variety of possible attractors in Section \ref{Results}. Furthermore, we identify conditions under which coupling can prevent essential extinction. We demonstrate two mechanisms by which the whole population can persist, whereas both subpopulations would undergo (essential) extinction without dispersal. Finally, we point out the special role of transients and crises in this model. We conclude with a discussion of the results in Section \ref{Conclusion}.
\section{Model} \label{Model}
We consider a spatially structured population model of a single species in discrete time. We assume that at each time step dispersal occurs after reproduction \citep{hastings1993complex, lloyd1995coupled}. The order of events, since there are only two, does not affect the dynamics.
\subsection{Reproduction (Local dynamics)}
The local dynamics are defined by the Ricker map \citep{ricker1954stock} combined with positive density dependence by an Allee effect. One way to model this is
\begin{linenomath}
\begin{align}
    f(x_t) = x_t e^{r\left(1-\frac{x_t}{K}\right)\left(\frac{x_t}{A}-1\right)} \label{f} \text{ ,}
\end{align}
\end{linenomath}
where $x_t$ is the population density at time step $t$ and $f(x_t)$ is the population production. Parameters $r$, $K$ and $A$ describe the intrinsic per-capita growth, the carrying capacity and the Allee threshold, respectively, $r>0$ and $0<A<K$.\\
Applications of this model can be found, for instance, in fisheries or insect models \citep{walters1976adaptive, turchin1990rarity, estay2014role}. While this model is not intended to be a realistic representation of a particular species \citep{neubert1997simple}, it captures the main biological features of interest, i.e. the Allee effect and overcompensation. As such, our model formulation, similar to \citet{schreiber2003allee}, satisfies the following properties:
\begin{itemize}
    \item There is a unique positive density $D$ that leads to the maximum population density $M$ in the next generation
    \item Extremely large population densities lead to extremely small population densities in the next generation
    \item Populations under the Allee threshold $A$ will go extinct
\end{itemize}
These conditions also hold for other models of that type, e.g. the logistic map with Allee effect or a harvesting term.\\
Our form of $f$ is chosen in such a way that the Allee threshold is at a fixed value. Other formulations which are based on biological mechanisms \citep{schreiber2003allee,courchamp2008allee} may be more realistic but make visualization more difficult. However, our results do not depend on this choice. 
\subsection{Dispersal (Between-patch dynamics)}
\newcommand{\R}{\mathbb{R}}
\newcommand{\C}{C}

We consider two patches with population densities $x_t$ and $y_t$ at time $t$. In each patch, we assume the same reproduction dynamics as in Equation (\ref{f}). The patches are linked by dispersal:
\begin{linenomath}
\begin{equation}
\begin{aligned}\label{xy}
    x_{t+1} &= (1-d) f(x_t) + d f(y_t) \text{ ,}\\
    y_{t+1} &= (1-d) f(y_t) + d f(x_t) \text{ ,}
\end{aligned}
\end{equation}
\end{linenomath}
where $d\in[0,0.5]$ is the fraction of dispersers ($0.5$ corresponds to complete mixing). Note that apart from initial conditions, the two patches are identical. The state space for this two-patch system is the non-negative cone $\C=[0,\infty)^2$ of $\R^2$. The solutions of \eqref{xy} correspond to iterating the map $F:\C\to\C$ given by $F(x,y)=((1-d)f(x)+df(y),df(x)+(1-d)f(y))$.

\section{Results} \label{Results}
\subsection{Dynamics without dispersal} \label{Local}
\begin{figure}
\centering
\begin{subfigure}{0.45\textwidth}
\includegraphics[width=\linewidth]{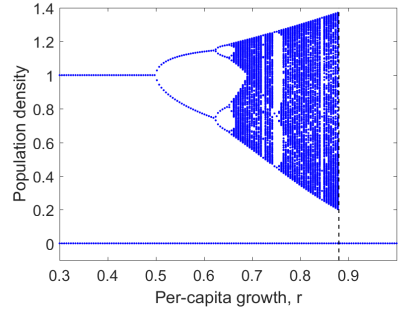}
\subcaption{Local dynamics}
\end{subfigure}\\
\medskip
\begin{subfigure}{0.45\textwidth}
\includegraphics[width=\linewidth]{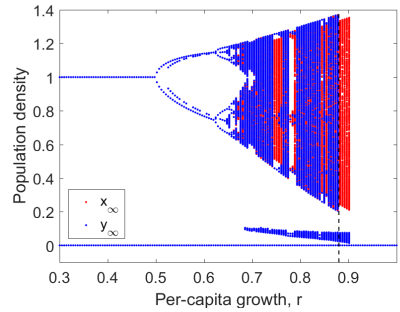}
\subcaption{Coupled system $d=0.03$}
\end{subfigure}\\
\medskip
\begin{subfigure}{0.45\textwidth}
\includegraphics[width=\linewidth]{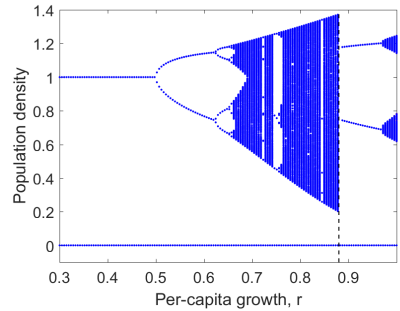}
\subcaption{Coupled system $d = 0.24$}
\end{subfigure}
\caption{Bifurcation diagram with bifurcation parameter $r$ of (a) the dynamics of a single isolated population and of two populations in the coupled system with (b) dispersal fraction $d=0.03$ whereby $x_\infty$ (red) is hidden partially by $y_\infty$ (blue) and (c) $d=0.24$ whereby $x_\infty = y_\infty$, thus only one patch is visible. The essential extinction threshold of an isolated population is marked with a dashed vertical line at $r_{th} = 0.88$. Dispersal can prevent extinction for $r > r_{th}$ in (b) and (c). Allee threshold $A = 0.2$, carrying capacity $K = 1$ and 8000 time steps of which the last 300 are plotted. Initial conditions: $(0.08,0.19)$, $(0.44,0.14)$, $(0.73,0.11)$, $(0.76,0.73)$, $(0.99,0.17)$ in all simulations.}

\label{fig:bif}
\end{figure}
In this section, we recap results from the local dynamics which are qualitatively similar to \cite{schreiber2003allee}. System (\ref{f}) has three equilibria, $x^*_1 = 0$, $x^*_2 = A$ and $x^*_3 = K$. We distinguish two dynamical patterns for the local case, depending on the threshold value $r_{th}$ that fulfills the equation $f(f(D)) = A$. For $0 < r < r_{th}$ the system is bistable. There is an upper bound $\bar{A}$ with $f(\bar{A})=A$. For initial densities $A<x_0<\bar{A}$, the population persists and goes extinct otherwise. The extinction attractor $x^*_1$ is always stable whereas the persistence attractor can be:
\begin{itemize}
\item A fixed point/an equilibrium for which $x_t = f(x_t)$;
\item A periodic orbit\footnote{Note that a fixed point is a periodic orbit of period one.} for which $x_t = f^n(x_t)$ but $x_t \neq f^j(x_t) \text{ } \forall \text{ } j=1,...,n-1$; or
\item A chaotic attractor \citep[see][for a definition]{broer2010dynamical}.
\end{itemize}
It loses its stability when $r>r_{th}$ and almost every initial density leads to essential extinction, i.e. for a randomly chosen initial condition with respect to a continuous distribution, extinction occurs with probability one \citep{schreiber2003allee}. This is shown in a bifurcation diagram with respect to $r$ in Figure \ref{fig:bif}a. The threshold $r_{th}$ is marked with a dashed line. These properties of the local dynamics (\ref{f}) can be formalized in a Theorem (see Appendix A).\\
Before turning towards the coupled model we consider two isolated patches, that is System (\ref{xy}) and $d=0$. For relatively small values of $r$ the persistence attractor of $f$ is a fixed point. The combination of equilibria of System (\ref{f}) delivers the equilibria of the uncoupled System (\ref{xy}): $(0,0)$, $(K,0)$, $(0,K)$, $(K,K)$, $(A,0)$, $(0,A)$, $(A,A)$, $(K,A)$ and $(A,K)$. Similar to \cite{amarasekare1998allee}, the last five equilibria are unstable. The first four equilibria are stable.\\
However,for larger values of $r$, the persistence attractor is not necessarily a fixed point and can be periodic or chaotic. When it has a linearly stable periodic orbit $\{p, f(p),\ldots , f^{n-1}(p)\}$ of period $n \geq 1$, the uncoupled map has $n + 3$ stable periodic orbits given by the forward orbits of the following periodic points
\begin{linenomath}
\begin{align}
\mathcal{P} = \{(0, 0),(0, p),(p, 0),(p, p),(p, f(p)), \ldots ,(p, f^{n-1}(p))\} \text{ .} \label{eq3}
\end{align}
\end{linenomath}
For the biological interpretation of the model it is important to note that one can obtain either global extinction of the whole population or persistence above the Allee threshold in one or both patches in the long term. The outcome follows from the dynamical behaviour of the local system. That changes with the introduction of dispersal. Attractors can appear or disappear and the fact that essential extinction always occurs for $r>r_{th}$ is no longer true.

\subsection{Additional attractors in the coupled system} \label{Multi}
\begin{figure}
\begin{subfigure}{0.48\textwidth}
\includegraphics[width=\linewidth]{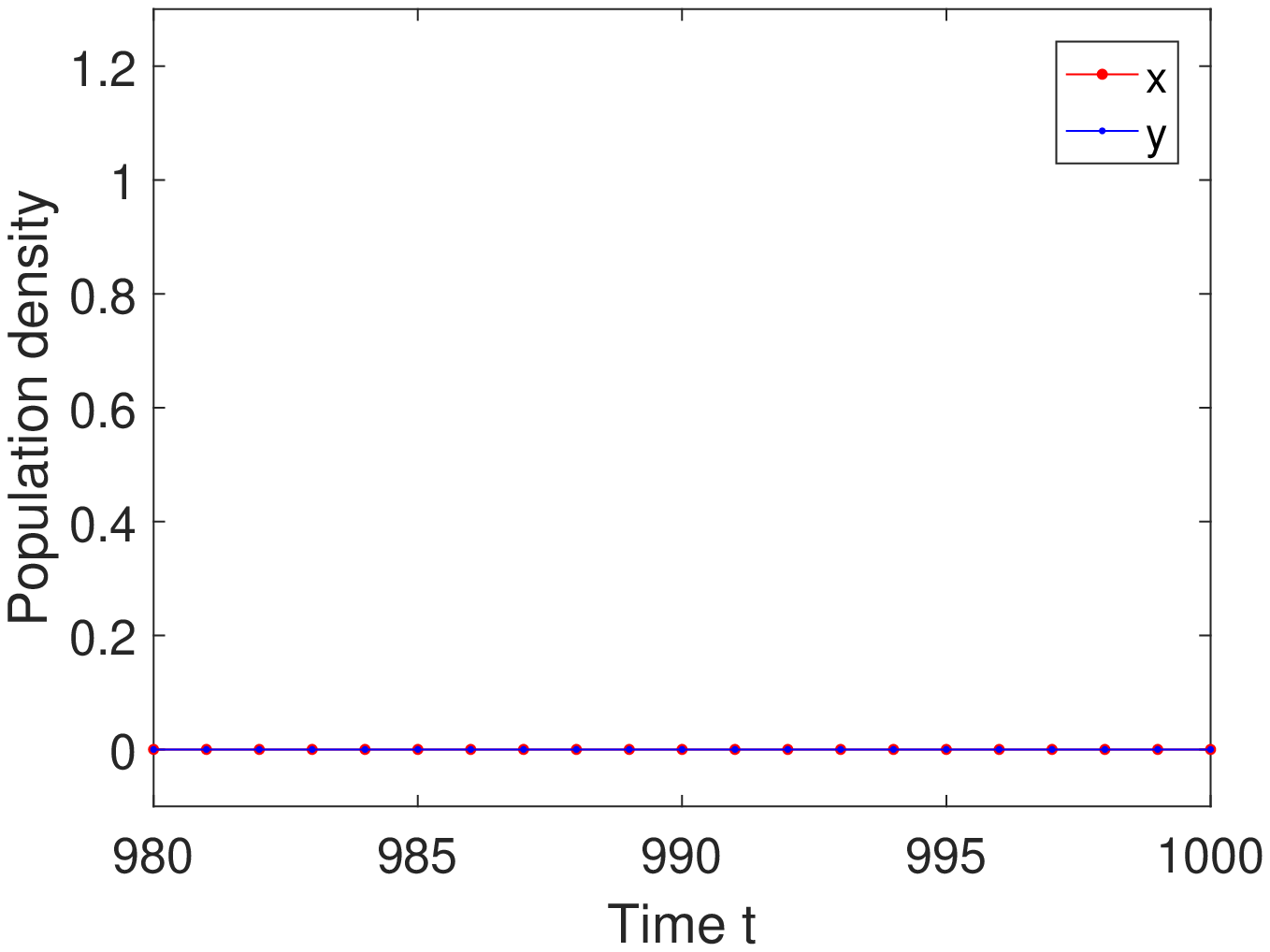}
\subcaption{Extinction}
\end{subfigure}\hspace*{\fill}
\begin{subfigure}{0.48\textwidth}
\includegraphics[width=\linewidth]{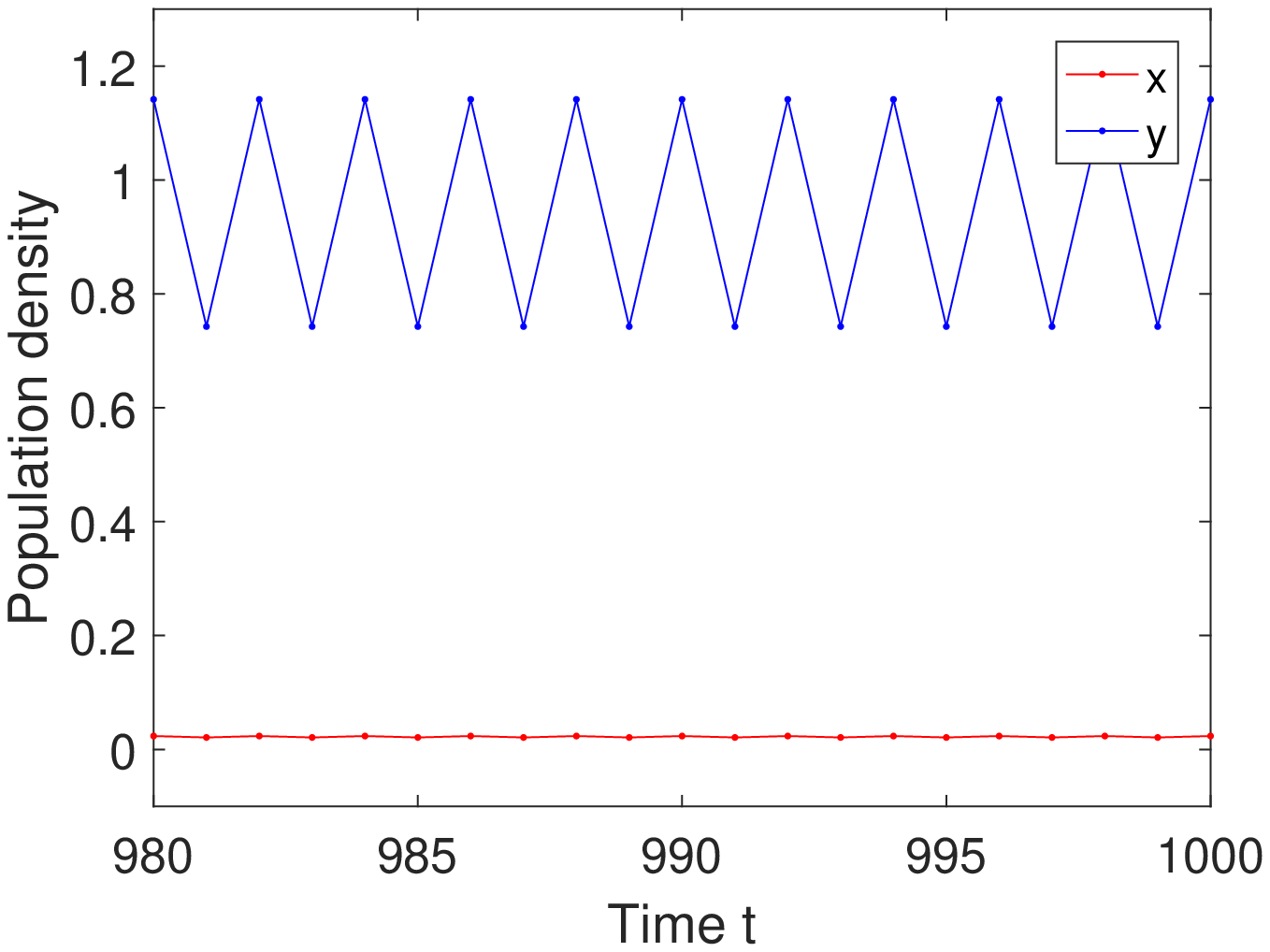}
\subcaption{Spatial asymmetry $y>x$}
\end{subfigure}

\medskip
\begin{subfigure}{0.48\textwidth}
\includegraphics[width=\linewidth]{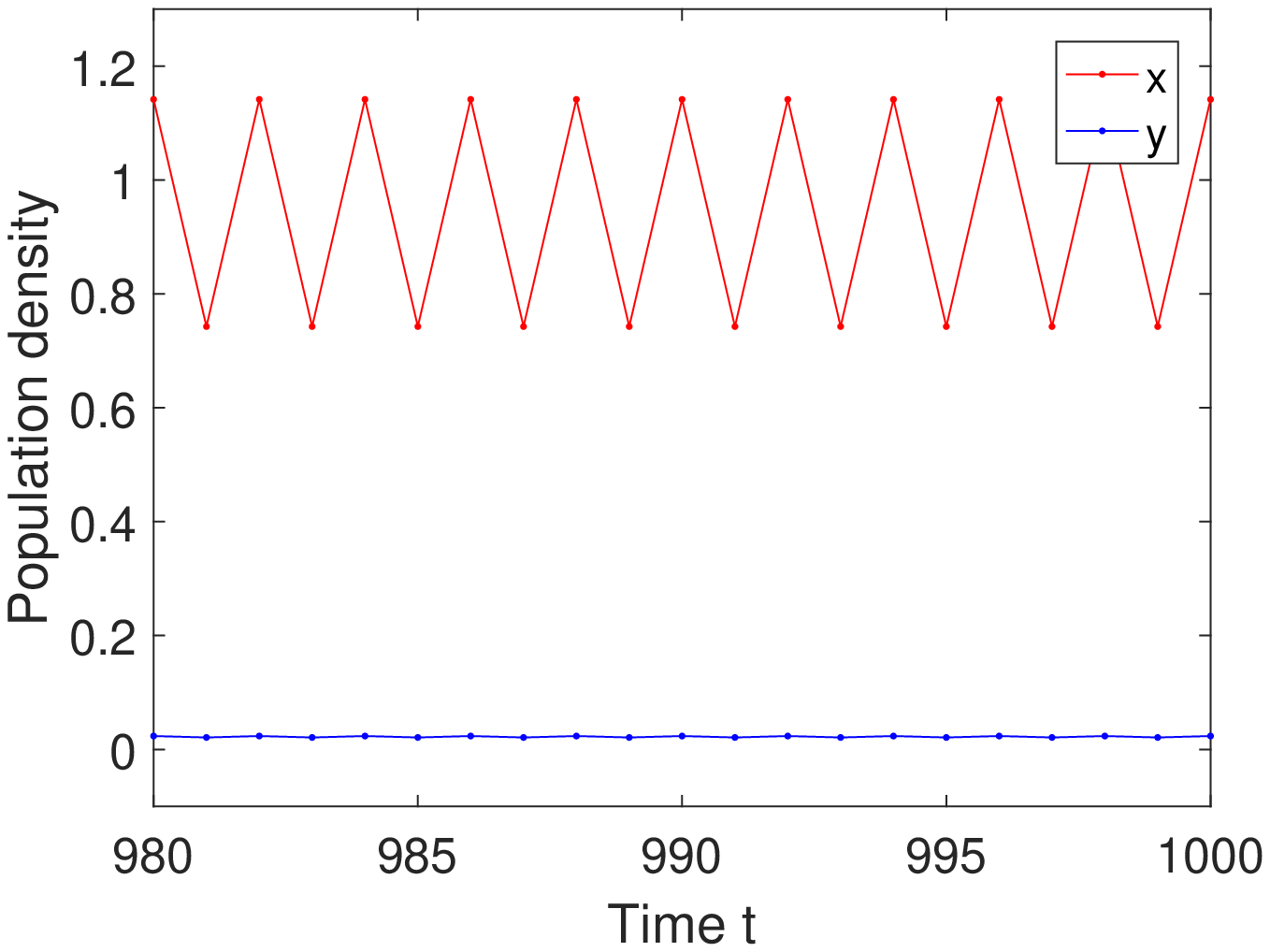}
\subcaption{Spatial asymmetry $x>y$}
\end{subfigure}\hspace*{\fill}
\begin{subfigure}{0.48\textwidth}
\includegraphics[width=\linewidth]{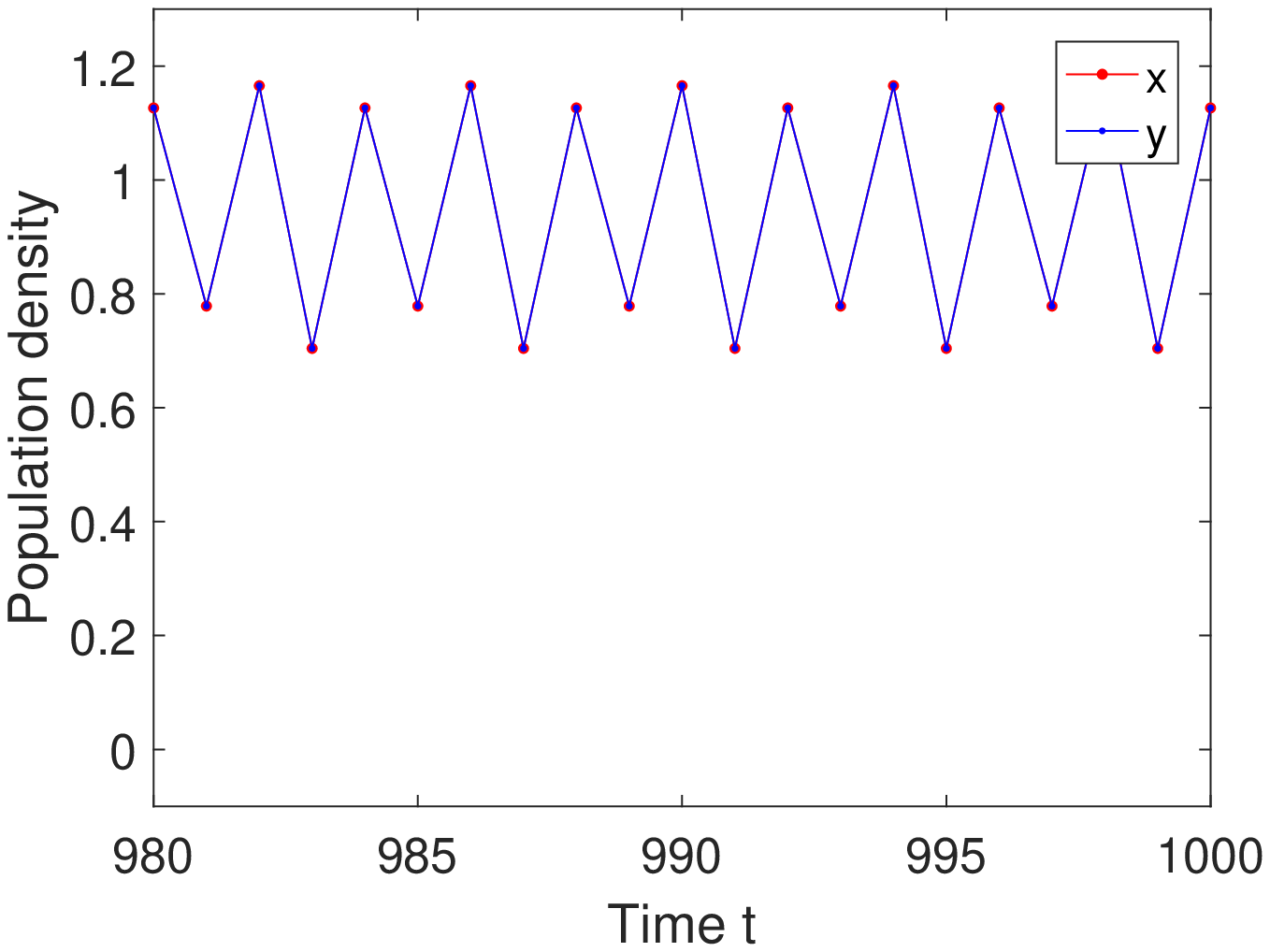}
\subcaption{In-phase 4-cycle}
\end{subfigure}

\medskip
\begin{subfigure}{0.48\textwidth}
\includegraphics[width=\linewidth]{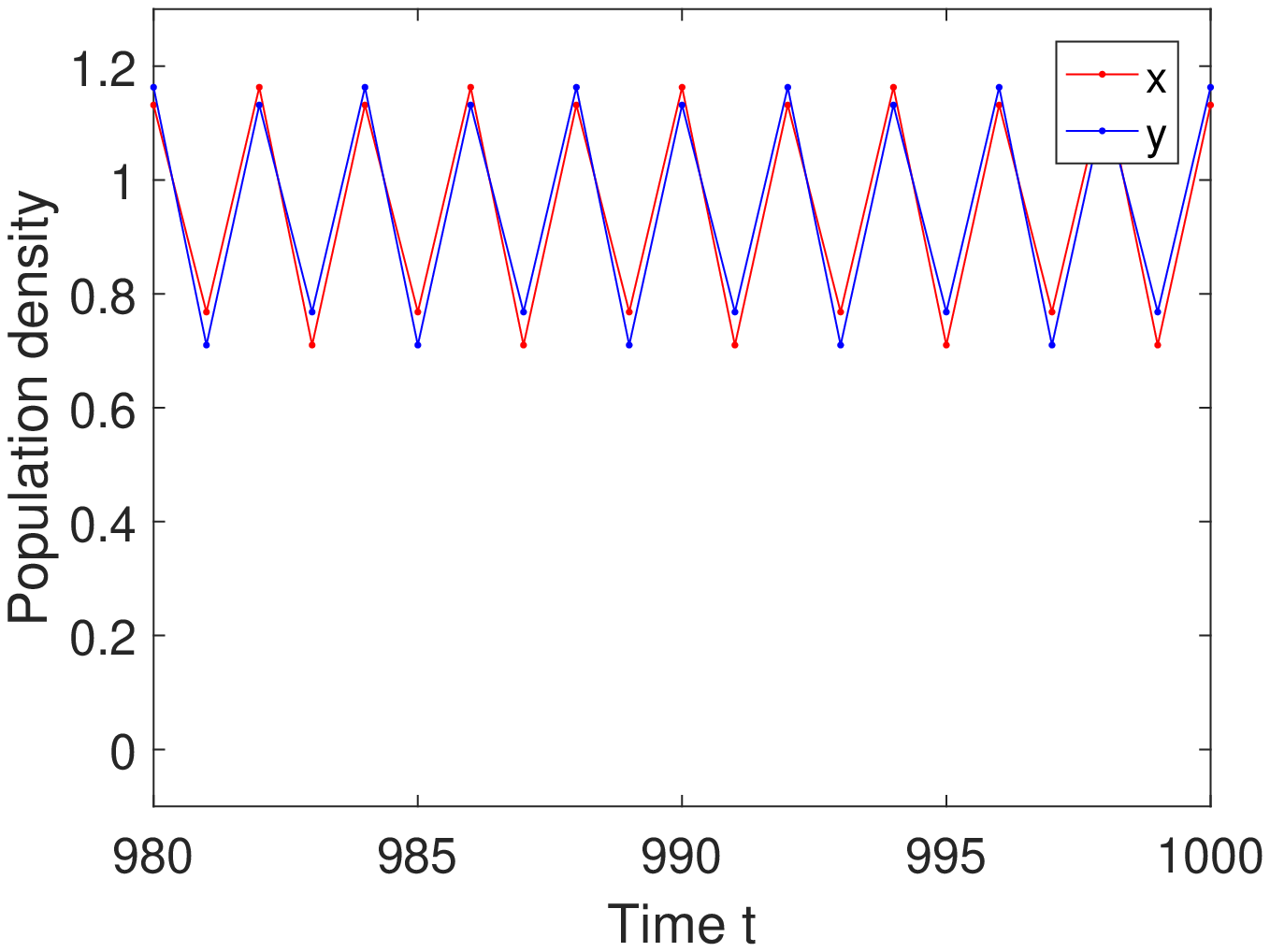}
\subcaption{Out-of-phase 4-cycle} 
\end{subfigure}\hspace*{\fill}
\begin{subfigure}{0.48\textwidth}
\includegraphics[width=\linewidth]{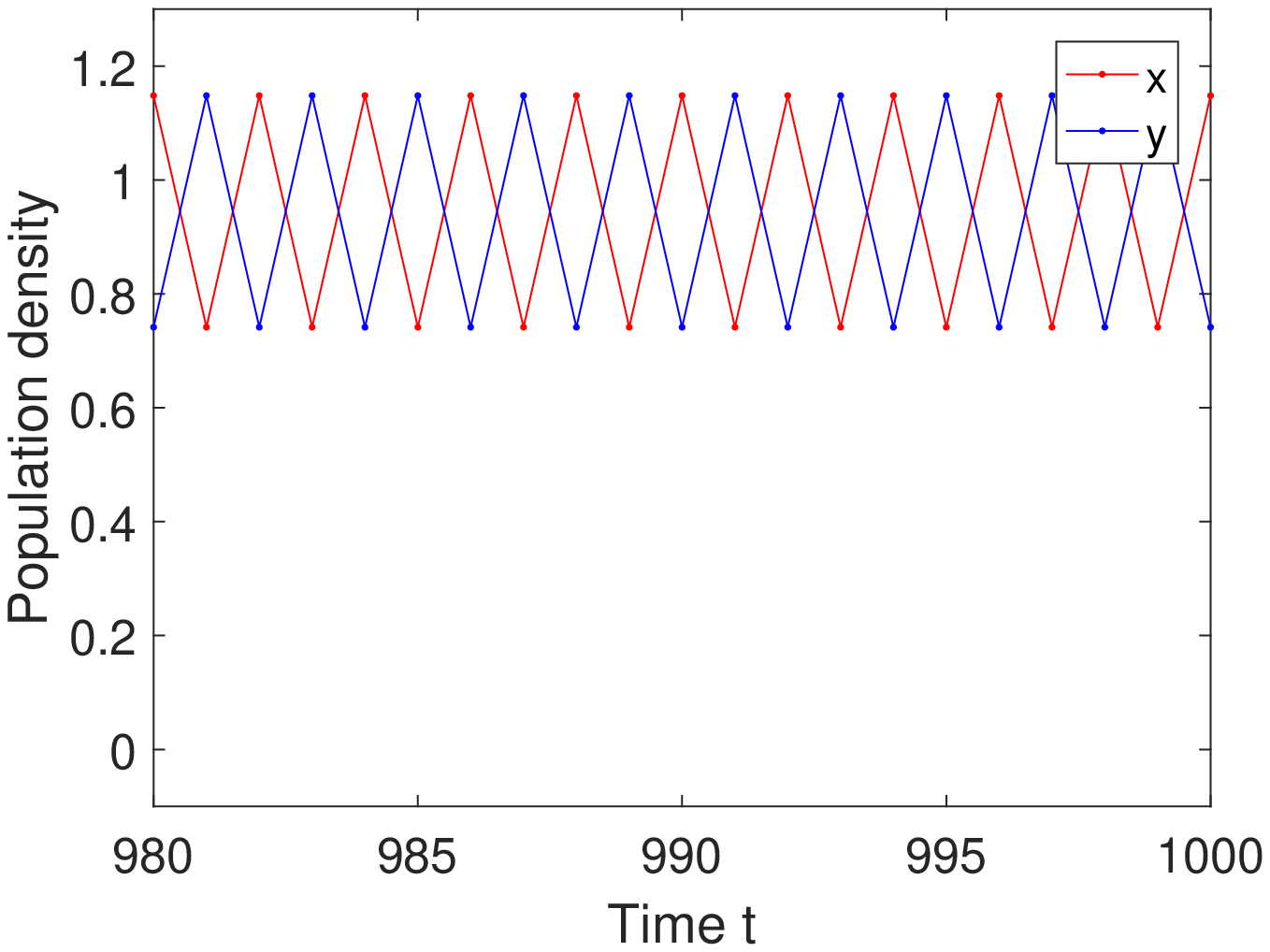}
\subcaption{Out-of-phase 2-cycle}

\end{subfigure}
\caption{Time series of model (\ref{xy}) that lead to different attractors because of different initial conditions. Parameters: $K=1$, $A=0.2$, $r=0.63$ and $d = 0.01$. Initial conditions: (a) $x_0=0.03$, $y_0=0.04$, (b) $x_0=0.16$, $y_0=0.86$, (c) $x_0=0.86$, $y_0=0.16$, (d) $x_0=0.64$, $y_0=0.38$, (e) $x_0=0.82$, $y_0=0.98$, (f)  $x_0=0.38$, $y_0=0.58$.} \label{fig:timeseries}
\end{figure}
When dispersal is weak and there is a stable positive periodic orbit for $f$, we prove the following theorem that shows that almost every initial condition converges to one of the $n+3$ stable periodic orbits in $\mathcal{P}$. Furthermore, if the positive stable periodic orbit of $f$ is not a power of $2$, then there are an infinite number of unstable periodic orbits.
\begin{theorem}\label{th2}
Assume the one-dimensional map $f(x)$ has a positive, linearly stable periodic orbit, $\{p, f(p), \ldots , f^{n-1}(p)\}$, with period $n \geq 1$. Let $U$ be an open neighborhood of $\cup_{i=1}^n (f \times f)^i(\mathcal{P})$. Then for $d > 0$ sufficiently small
\begin{itemize}
\item[(i)] System (\ref{xy}) has $n + 3$ distinct, linearly stable periodic orbits contained in
$U$. Let $G$ denote the union of these linearly stable periodic orbits.
\item[(ii)] $\C\setminus B$ has Lebesgue measure zero where\\
$B=\{(x,y)\in\C: \lim_{t\to\infty}\rm{dist}$$(F^t(x,y),G)=0$$\}$ is the basin of attraction of $G$.
\item[(iii)] If $n$ is not a power of $2$, then $\C\setminus B$ contains an infinite number of periodic points.
\end{itemize}
\end{theorem}
A proof of this theorem is given in Appendix B. Since $f$ is known to undergo period doublings until chaos, one can obtain a large number of attractors for weakly coupled maps. However, our numerical results show that for larger $d>0$, the number of coexisting attractors is smaller than $n+3$.\\ 
Consider System (\ref{xy}) with parameter values $r=0.63$ and $d=0.01$. This value of $r$ leads to 4-cycles in the uncoupled system. We observe six stable periodic orbits. Time series for different initial conditions are shown in Figure \ref{fig:timeseries}. The extinction state in both patches is stable (Figure \ref{fig:timeseries}a). The two attractors in Figure \ref{fig:timeseries}b and Figure \ref{fig:timeseries}c show periodic behaviour above the Allee threshold in one patch and below the Allee threshold in the other patch. We call these attractors \emph{asymmetric attractors}. \\
In contrast to four different 4-cycles for sufficiently small $d$ (Theorem \ref{th2}), we observe an in-phase 4-cycle (Figure \ref{fig:timeseries}d) and only one out-of-phase 4-cycle (Figure \ref{fig:timeseries}e). The other two 4-cycles
with $x_t<1$, $y_t>1$ and $x_{t+1}>1$, $y_{t+1}<1$ are replaced by only one attractor, an out-of-phase 2-cycle (Figure \ref{fig:timeseries}f). This is an example for a stabilizing effect of dispersal. In the following, we will call all attractors with population densities above the Allee threshold in both patches \emph{symmetric attractors}.\\
Final-state sensitivity depending on the initial conditions can occur whenever there are several coexisting attractors \citep{peitgen2006chaos}. The system can exhibit very different dynamic behaviours even if all parameter values are fixed \citep{lloyd1995coupled}. In the following sections, we will first categorize attractors in terms of subpopulations being above or below the Allee threshold. Secondly, we take a closer look at different symmetric attractors, like the ones in Figure \ref{fig:timeseries}d-f.\\

\begin{figure}[h!]
    \centering
  \includegraphics[width=0.9\linewidth]{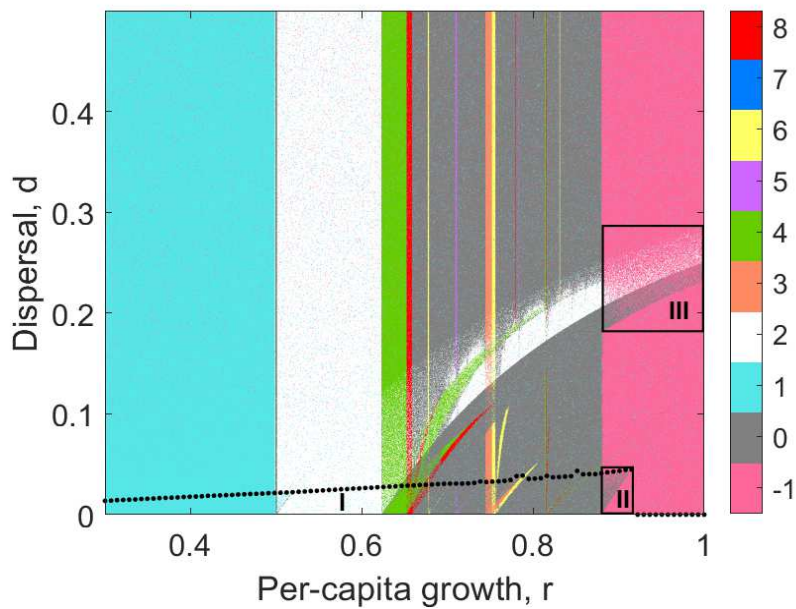}
\caption{Dynamical behaviour characterized by the periodicity, as a function of $r$ and $d$. Labels of the colour bar give the periodicity of locally stable cycles. Periodicity 1 stands for a stable equilibrium (trivial or non-trivial), 0 for periods $>8$ or chaos and -1 for extinction when $A < x_0 < 1 \vee A < y_0 < 1$. Region (I) below the dotted curve indicates for which values of $r$ and $d$ asymmetric attractors appear (tested for 100 random initial conditions) with irregularities due to additional attractors depending on dispersal. Regions (II) and (III) indicate for which values of $r$ and $d$ dispersal can prevent essential extinction. Fuzzy regions indicate multistability. Note that the extinction state is always stable (turquoise sprinkles). $K=1$ and $A=0.2$ fixed in all runs. One random initial condition per parameter combination. Selected periodicity has been determined using the \emph{CompDTIMe} routine for Matlab (\url{https://www.imath.kiev.ua/~nastyap/compdtime.html}), provided there was no essential extinction.}
\label{fig:app}
\end{figure}

For the simulations, we normalize the population density relative to the carrying capacity by setting $K=1$ and fix $A=0.2$. Then, there are only two remaining parameters, $r$ and $d$. Figure \ref{fig:app} summarizes the dynamical behaviour that can be observed in the ($r,d$)-parameter plane for $0 < d < 0.5$ and $0.3<r<1$. 
\subsubsection{Multiple attractors due to the Allee effect} \label{Spatial_Hetero}
In the case of weak dispersal (see Figure \ref{fig:app}, parameter region I, below dotted curve), the equilibria of the coupled system are similar to the ones of the uncoupled system. This follows from a perturbation argument, similar to \cite{karlin1972polymorphisms}. We observe four attractors that differ in whether the population density in each patch is above or below the Allee threshold. The extinction state $(0,0)$ is always stable. The two asymmetric and the symmetric attractors can be either equilibria or show periodic/chaotic behaviour, depending on the values of $r$ and $d$ (see Figure \ref{fig:app}). Thus, spatial asymmetry can be conserved. Figure \ref{fig:bif}b shows the four states in patch $y$ for $d=0.03$ (blue): when both subpopulations start above the Allee threshold, the population densities remain at carrying capacity $K$ or after period doublings on a periodic/chaotic attractor. If the initial population in patch $y$ is smaller than $A$ but larger in patch $x$, one asymmetric attractor is approached (red: large $x$, blue: small $y$). If initial populations in both patches are smaller than $A$, the extinction attractor is approached.\\
The situation changes for larger dispersal (see Figure \ref{fig:app}, above dotted curve). The asymmetric attractors disappear and only extinction or persistence above the Allee threshold in both patches is possible. This can be seen in Figure \ref{fig:bif}c, where in comparison to Figure \ref{fig:bif}b no asymmetric attractor is visible.\\
A nullcline analysis can give information about the number of equilibria that can lead to different attractors. For that, we refer to \cite{amarasekare1998allee} or \cite{kang2011expansion}, who did a detailed nullcline analysis for a corresponding continuous-time model. 

\subsubsection{Multiple attractors due to overcompensation} \label{Attractors2}
Multiple attractors can not only appear due to Allee effects but also in coupled maps with overcompensation \citep{hastings1993complex}. Thus, we take a closer look at additional symmetric attractors as shown in Figure \ref{fig:timeseries}d-f. The in-phase 4-cycle, the out-of-phase 4-cycle and the out-of-phase 2-cycle can coexist even without additional equilibria.\\
The ($r,d$)-parameter plane in Figure \ref{fig:app} provides some insights for which parameter combinations multiple symmetric attractors appear (Note that in this Figure, we do not distinguish between different attractors of the same period for better clarity): On the one hand, the equilibrium $(K,K)$ undergoes several period-doublings up to chaos and finally essential extinction when increasing $r$, independently of dispersal (vertical stripe structure). The bending stripes across the diagram, on the other hand, indicate additional attractors depending on both $r$ and $d$. Fuzzy regions appear when multiple symmetric attractors coexist. Coexisting symmetric attractors can be also seen in Figure \ref{fig:bif}b for $0.5 < r < 0.65$ where in-phase and out-of-phase 2-cycles coexist.\\
This phenomenon is well understood in models without an Allee effect \citep{hastings1993complex, yakubu2002interplay, wysham2008sudden, yakubu2008asynchronous}. As it only occurs for the symmetric attractor, where we observe population densities above the Allee threshold, the Allee effect itself is negligible concerning the origins of the non-equilibrium attractors. However, it is important to mention here, since any of the coexisting attractors can disappear due to the Allee effect with the system then collapsing to the extinction attractor. This will be discussed in Sections \ref{DIPEE} and \ref{Transients}.\\
Combining the results of discrete-time models with overcompensation \citep{hastings1993complex, lloyd1995coupled, kendall1998spatial} and continuous-time models for spatially structured populations with Allee effect \citep{amarasekare1998allee} shows that the variety of both is expressed here.

\subsection{Dispersal induced prevention of essential extinction} \label{DIPEE}

In Section \ref{Local} we have seen that for per-capita growth exceeding the threshold $r_{th}$ isolated populations undergo essential extinction. We now investigate mechanisms that allow ``dispersal induced prevention of essential extinction" (DIPEE) in the coupled maps. We choose the parameters such that without dispersal the whole population would go extinct ($r > r_{th}$). We identify two mechanisms for DIPEE: Spatial asymmetry and stabilizing (approximately) out-of-phase dynamics.
\subsubsection{DIPEE due to spatial asymmetry}
\begin{figure}
\begin{subfigure}{0.43\textwidth}
\includegraphics[width=\linewidth]{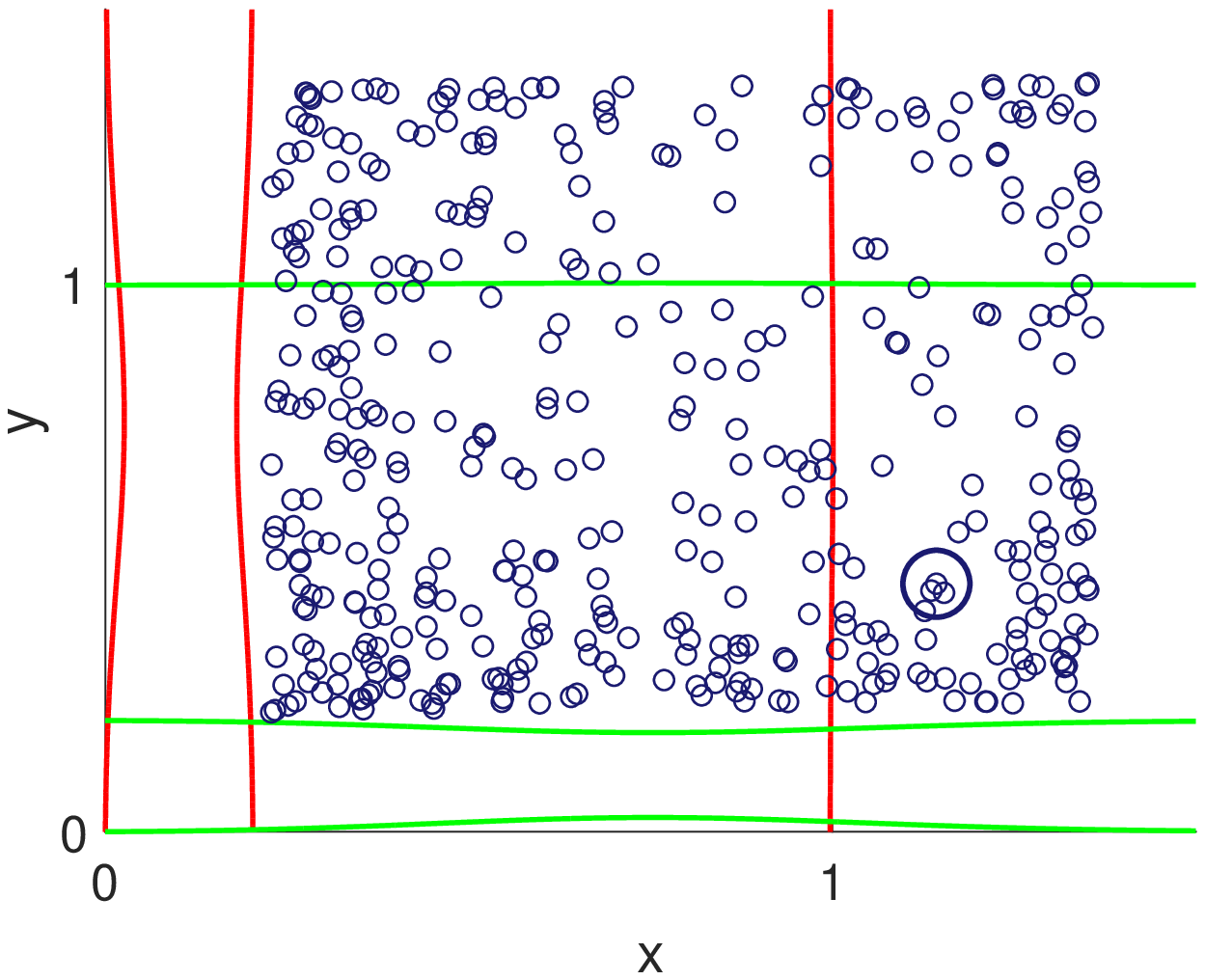}
\subcaption{$r=0.87 < r_{th}$}
\end{subfigure}\hspace*{\fill}
\begin{subfigure}{0.43\textwidth}
\includegraphics[width=\linewidth]{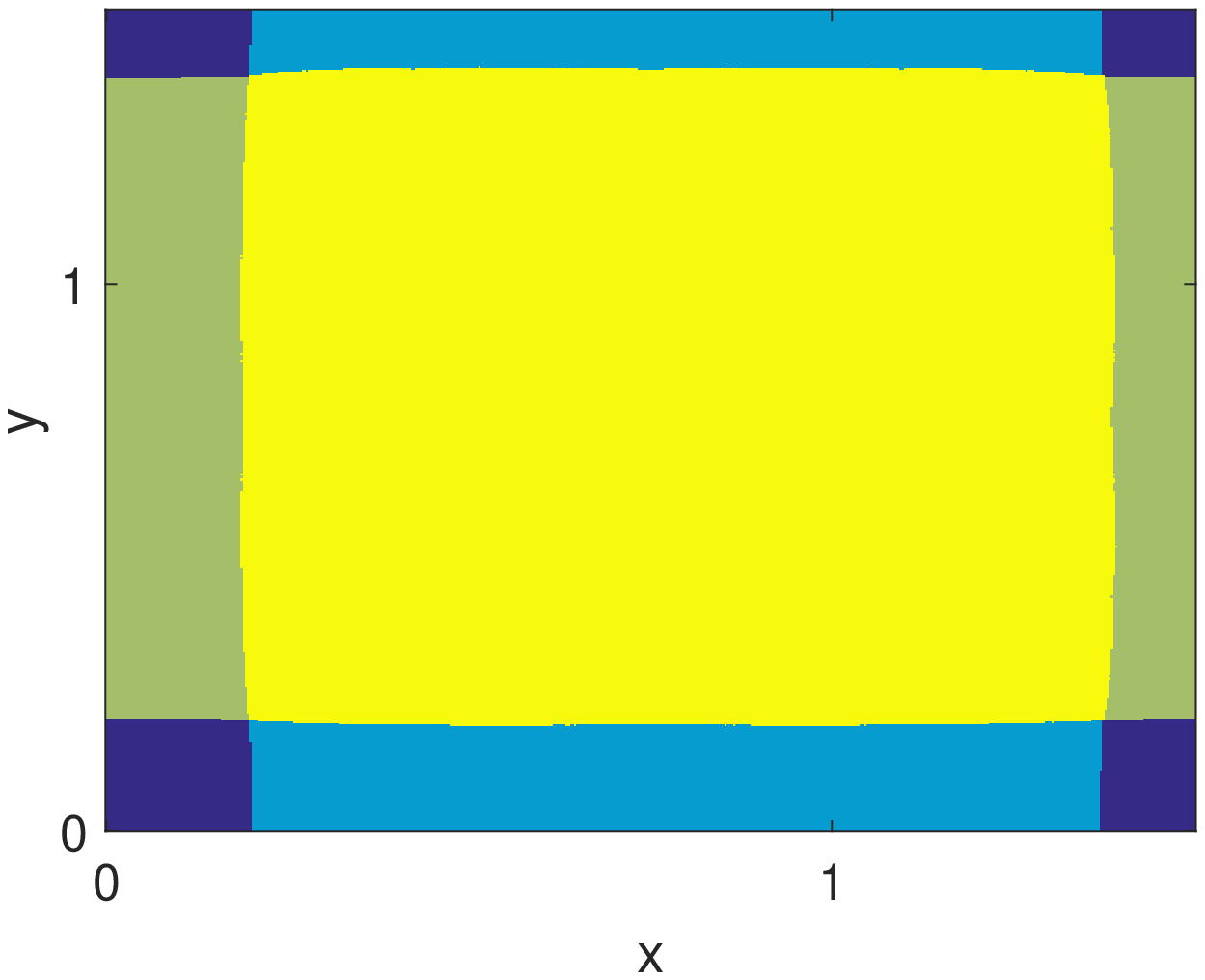}
\subcaption{$r=0.87 < r_{th}$}
\end{subfigure}

\medskip
\begin{subfigure}{0.43\textwidth}
\includegraphics[width=\linewidth]{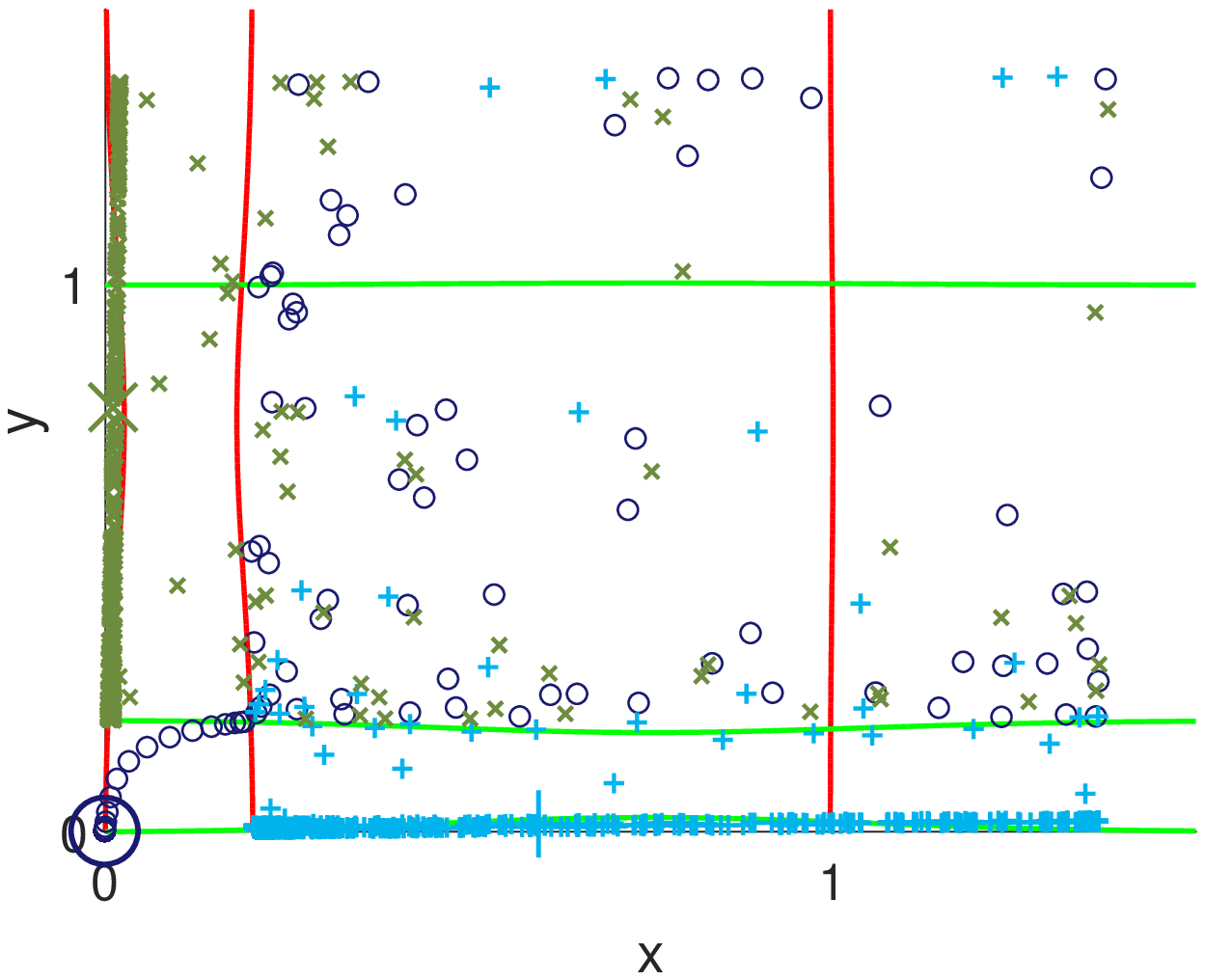}
\subcaption{$r=0.887 > r_{th}$}
\end{subfigure}\hspace*{\fill}
\begin{subfigure}{0.43\textwidth}
\includegraphics[width=\linewidth]{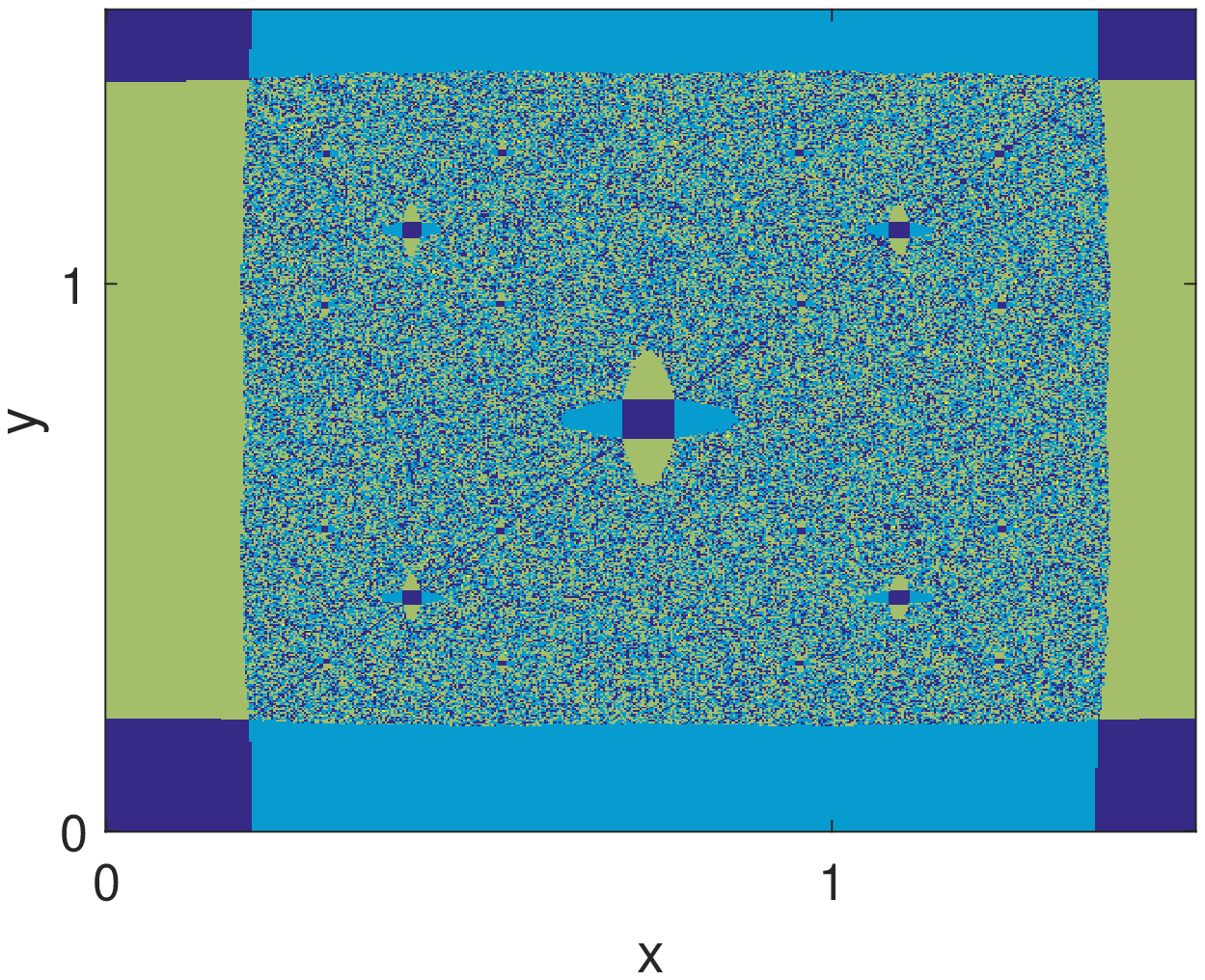}
\textbf{\subcaption{$r=0.887 > r_{th}$}}
\end{subfigure}

\medskip
\begin{subfigure}{0.43\textwidth}
\includegraphics[width=\linewidth]{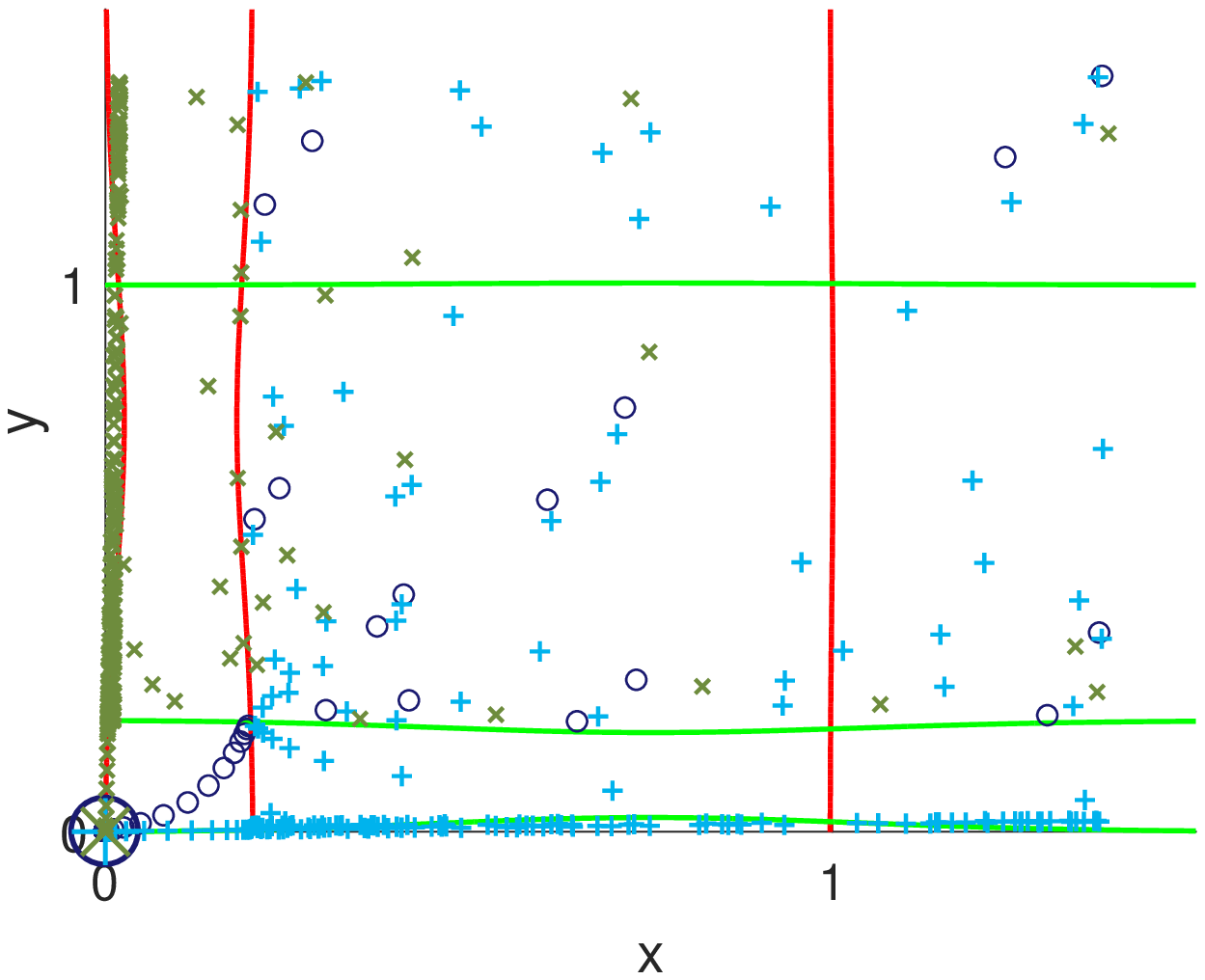}
\subcaption{$r=0.888 > r_{th}$}
\end{subfigure}\hspace*{\fill}
\begin{subfigure}{0.43\textwidth}
\includegraphics[width=\linewidth]{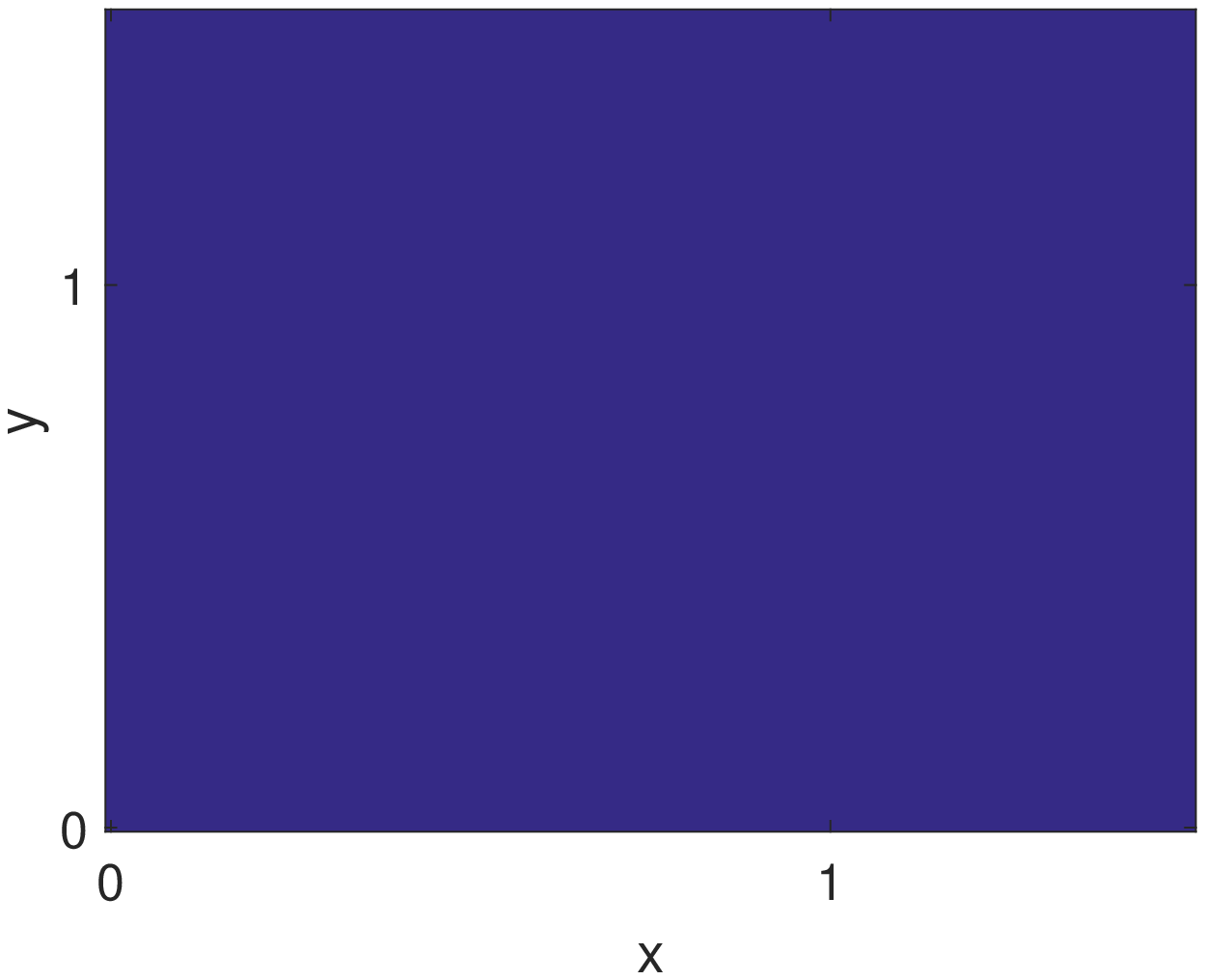}
\subcaption{$r=0.888 > r_{th}$}
\end{subfigure}
\caption{Phase planes (left column) and basins of attraction (right column) of the coupled system with $d=0.01$ and (a,b) $r=0.87$, (c,d) $r=0.887$ and (e,f) $r=0.888$. In the phase planes, sample orbits for initial conditions $(A,A) < (x_0,y_0) < (\bar{A},\bar{A})$ are shown with dots/crosses. When $r < r_{th}$ the population persists (a). For $r$ exceeding $r_{th}$ two asymmetric states (and thus DIPEE) and the extinction state are possible (c). For sufficiently large $r$ extinction is inevitable (e). Large symbols mark the final states. Nullclines in red and green, respectively. Basins of attraction of the four attractors in panels (b), (d) and (f): Extinction (dark blue), asymmetric coexistence (light blue and green), symmetric coexistence (yellow). Clear basin boundaries (b), fractal basin boundaries between asymmetric coexistence and extinction (d) or no boundaries (f) depending on the value of $r$. Allee threshold $A = 0.2$, carrying capacity $K = 1$ and 2000 time steps in all simulations.} \label{fig:test2}
\end{figure}

For the moment, we only consider small dispersal $d<0.05$ (see Figure \ref{fig:app}, parameter region II). In this case, the coupling is sufficiently weak to observe different dynamics in both patches. Figure \ref{fig:test2}a, \ref{fig:test2}c and \ref{fig:test2}e show the phase planes with nullclines\footnote{The $x$-nullcline is the set of points satisfying $x_{t+1} = x_t$, cf. \cite{kaplan2012understanding}. Similarly, the $y$-nullcline satisfies $y_{t+1} = y_t$} and sample orbits for different values of $r$. In Figure \ref{fig:test2}a all orbits with initial conditions $(A,A) < (x_0,y_0) < (\bar{A},\bar{A})$ remain on the chaotic symmetric attractor. When $r$ exceeds $r_{th}$, the symmetric attractor collides with the unstable equilibrium $(A,A)$ and disappears whereas the asymmetric attractors persist. \cite{grebogi1982chaotic} and \cite{bischi2016qualitative} call that phenomenon a boundary crisis. Figure \ref{fig:test2}c presents three cases of orbits with $(A,A) < (x_0,y_0) < (\bar{A},\bar{A})$: either the whole population goes extinct (dark blue) or the population in one patch drops under the Allee threshold, while the population in the other patch remains above (light blue, green). In this situation, essential extinction can be prevented, depending on the initial conditions. One subpopulation overshoots the equilibrium beyond some critical value (e.g. in patch $x$) and then drops below the Allee threshold whereas the other subpopulation (e.g. patch $y$) remains above. This leads to high net dispersal from patch $y$ to patch $x$. Thus, in patch $y$, the maximum population density is reduced, so that $f(M)>A$ and essential extinction does not take place. Patch $x$ is rescued from extinction by continual migration from patch $y$.\\
The basins of attraction change when $r$ exceeds $r_{th}$. For $r<r_{th}$ the basins are sharply separated sets as shown in Figure \ref{fig:test2}b. When the symmetric attractor disappears its basin results in a fractal structure (see Figure \ref{fig:test2}d). When parameter $r$ is increased further, DIPEE is not possible. The two asymmetric attractors disappear after another boundary crisis with equilibria near $(0,A)$ and $(A,0)$ (see Figure \ref{fig:test2}e). Almost all initial conditions lead to the only remaining attractor, the extinction state (see Figure \ref{fig:test2}f).\\
In summary, for per-capita growth above the local essential extinction threshold $r_{th}$ small dispersal can have a stabilizing effect in terms of reducing the maximum population density and thus preventing essential extinction (see Figure \ref{fig:app}, parameter region II). This result is emphasized by Figure \ref{fig:bif}b. The asymmetric attractor in which patch $y$ remains below and $x$ above the Allee threshold can persist for values of $r>r_{th}$. Conversely, one can observe the symmetric attractor to disappear at $r_{th}$. Note that the opposite case in which patch $x$ is below $A$ also persists for $r>r_{th}$ but is not shown in Figure \ref{fig:bif}b.

\subsubsection{DIPEE due to stabilizing (approximately) out-of-phase dynamics} \label{Opposed}
\begin{figure}
\begin{subfigure}{0.5\textwidth}
\includegraphics[width=\linewidth]{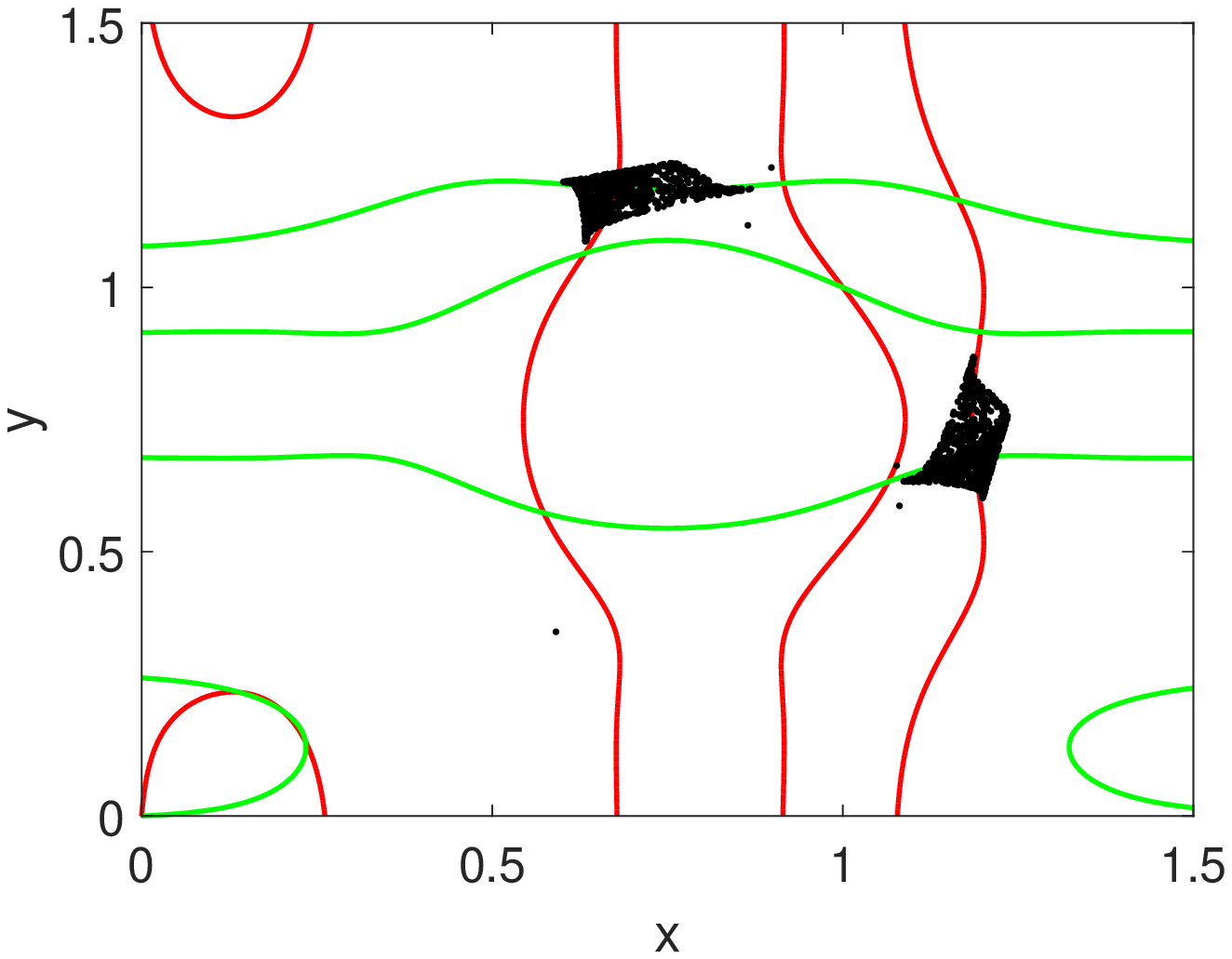}
\caption{$d=0.19$}
\end{subfigure}\hspace*{\fill}
\begin{subfigure}{0.5\textwidth}
\includegraphics[width=\linewidth]{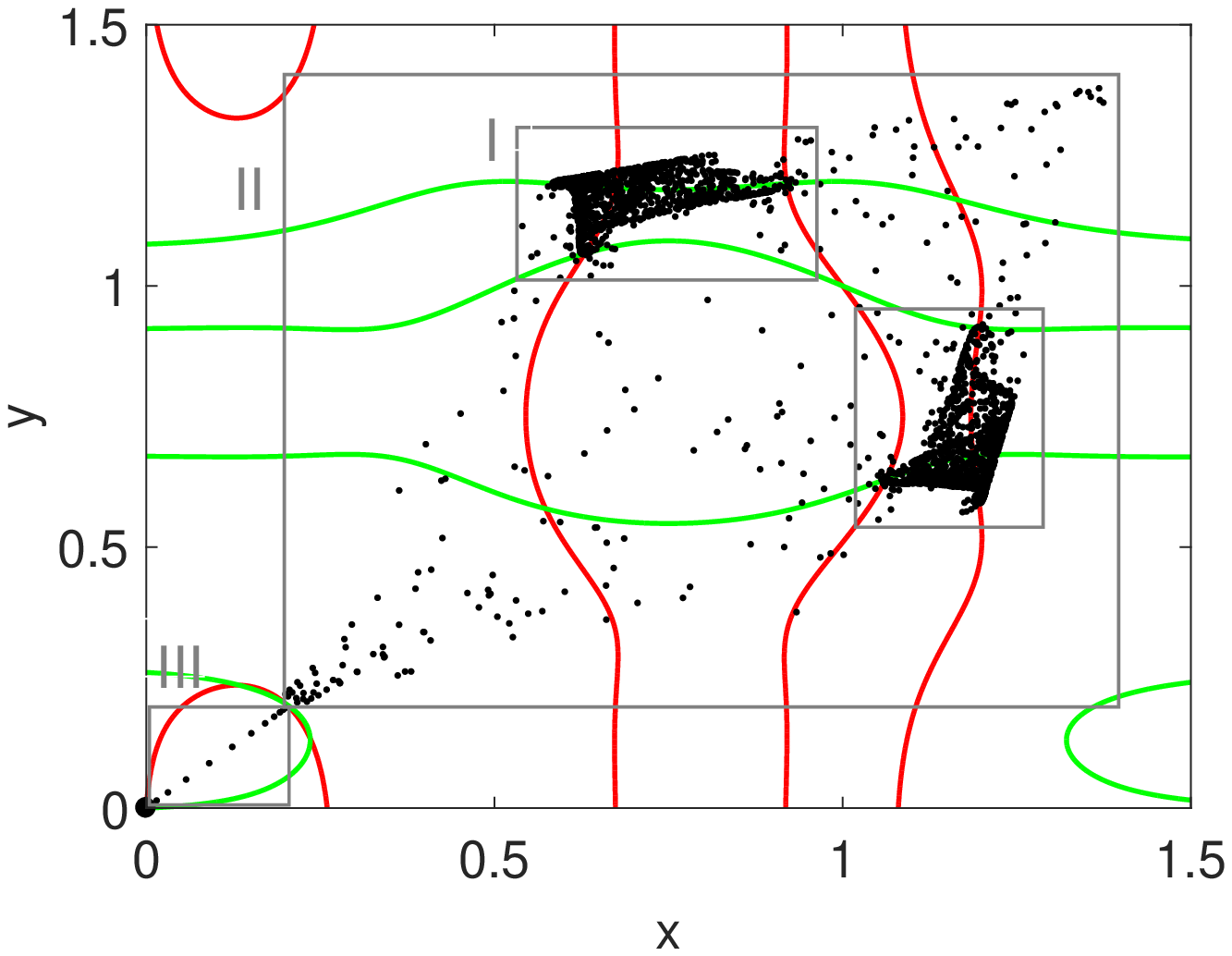}
\caption{$d=0.186$}
\end{subfigure}
\caption{Phase planes with nullclines of the second iteration of System (\ref{xy}) with $r=0.89$ and (a) $d=0.19$ and (b) $d=0.186$. The approximately out-of-phase attractor (a) undergoes a boundary crisis (b, region I). The emerging chaotic rhombus (b, region II) again merges the unstable equilibrium $(A,A)$ and finally converges to the extinction state (b, region III). Allee threshold $A = 0.2$, carrying capacity $K = 1$, $(A,A) < (x_0,y_0) < (\bar{A},\bar{A})$, 1000 time steps, large symbol: final state.} \label{fig:test5}
\end{figure}
\begin{figure}
\begin{subfigure}{0.5\textwidth}
\includegraphics[width=\linewidth]{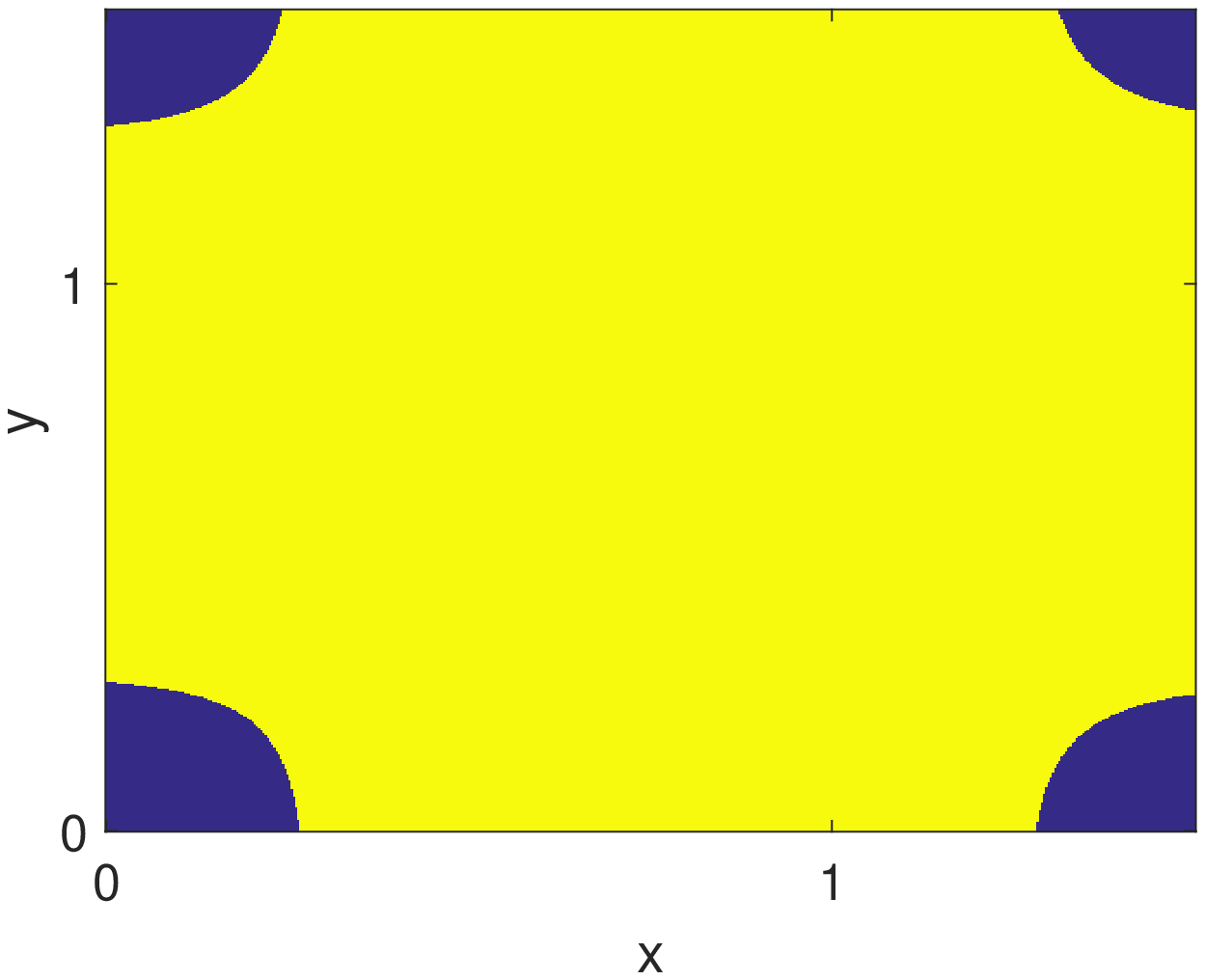}
\caption{$r = 0.87$}
\end{subfigure}\hspace*{\fill}
\begin{subfigure}{0.5\textwidth}
\includegraphics[width=\linewidth]{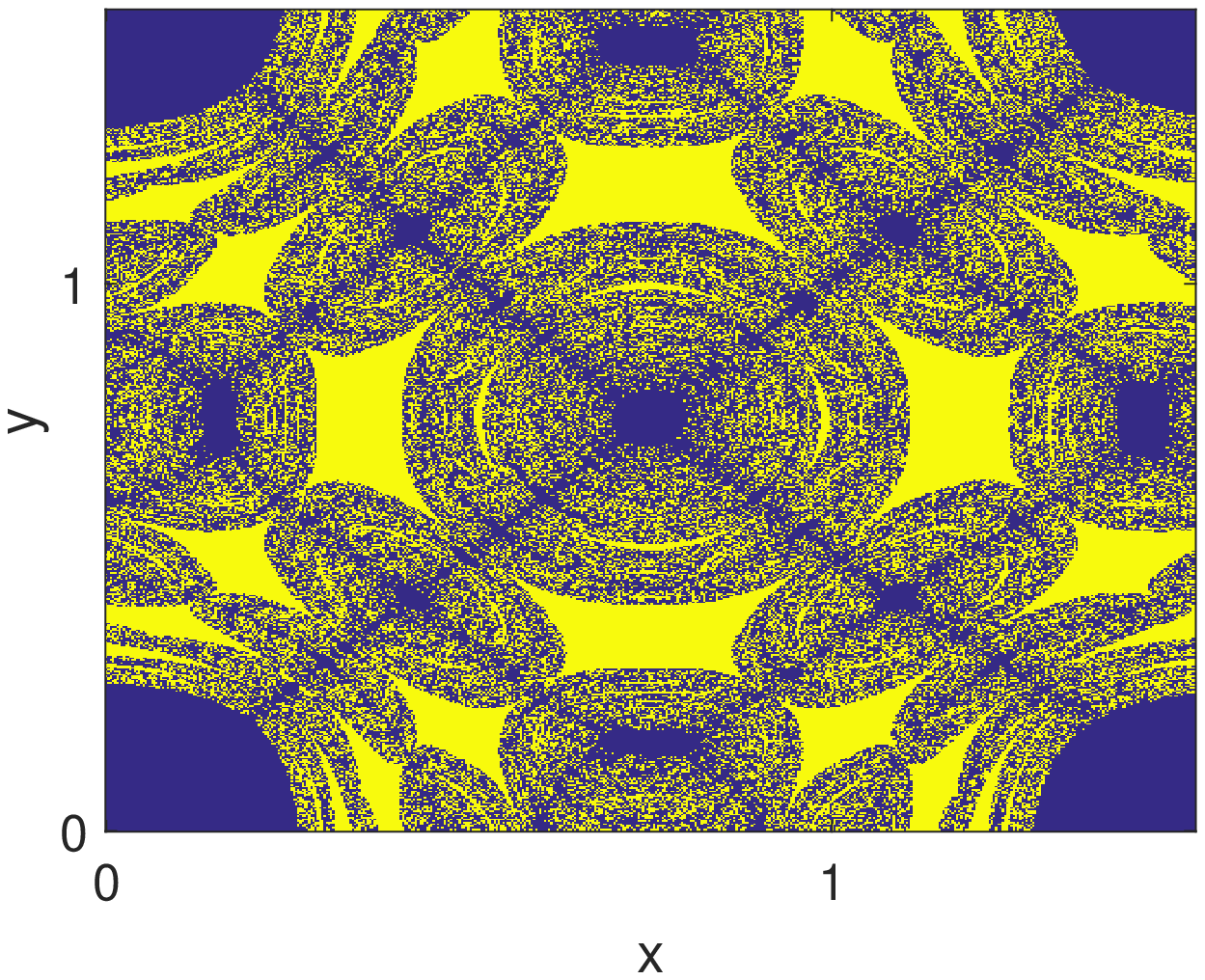}
\caption{$r=0.887$}
\end{subfigure}
\caption{Basins of attraction for $d=0.23$ and (a) $r=0.87$ and (b) $r=0.887$. Blue indicates the extinction state whereas yellow marks symmetric coexistence attractors. When $r$ exceeds $r_{th}$, the basins change to a fractal structure. Allee threshold $A = 0.2$, carrying capacity $K = 1$ and 1000 time steps in all simulations.} \label{fig:bas_attr2}
\end{figure}

A second mechanism that can prevent essential extinction operates at larger dispersal fractions around $0.19<d<0.28$ (see Figure \ref{fig:app}, parameter region III). In this parameter region, asymmetric attractors are impossible. Both subpopulations either persist above the Allee threshold or go extinct. The extinction state $(0,0)$ is stable whereas the symmetric attractor shows (approximately) out-of-phase dynamics where both population densities are above the Allee threshold but alternating (see Figure \ref{fig:test5}a). For values $r<r_{th}$, the symmetric out-of-phase dynamics coexist with a chaotic rhombus\footnote{no obvious relationship between $x_t$ and $y_t$ \citep{kendall1998spatial}; the attractor forms a rhombic structure}. Initial conditions $(A,A) < (x_0,y_0) < (\bar{A},\bar{A})$ lead either to one or the other attractor. When $r$ exceeds $r_{th}$, the chaotic rhombus collides with the unstable equilibrium $(A,A)$ (similar to Figure \ref{fig:test3}) and disappears whereas the (approximately) out-of-phase dynamics persists. Figure \ref{fig:bif}c shows the drastic change of possible attractors at $r_{th}$. In one time step more individuals move from patch $x$ to $y$. In the next step, net movement is from $y$ to $x$ so that values in the two patches cover the same range. Thus, only one patch is visible in Figure \ref{fig:bif}c. The other patch is overlaid completely. The antagonistic net movement prevents an overshoot in both patches and both are rescued from essential extinction. Again, one should note that DIPEE is very sensitive to the choice of initial conditions. More precisely, different initial conditions $(A,A) < (x_0,y_0) < (\bar{A},\bar{A})$ lead either to synchronization and thus essential extinction or to coexistence with population densities above the Allee threshold in both patches and thus DIPEE.\\
Also the basins of attraction change when $r$ exceeds $r_{th}$. For $r<r_{th}$ the basins are sharply separated sets as shown in Figure \ref{fig:bas_attr2}a. When the chaotic rhombus disappears the basins of attraction for symmetric attractors split into a fractal structure (see Figure \ref{fig:bas_attr2}b). This structure is well known from other studies on coupled maps with local overcompensation \citep{gyllenberg1993does, hastings1993complex, lloyd1995coupled}. The significant difference here is that attractors are distinguished not in their period but in the sense that slightly different initial conditions lead either to survival or to extinction. From the ecological point of view, that is a crucial difference.\\
\subsubsection{No DIPEE}
\begin{figure}
\begin{subfigure}{0.5\textwidth}
\includegraphics[width=\linewidth]{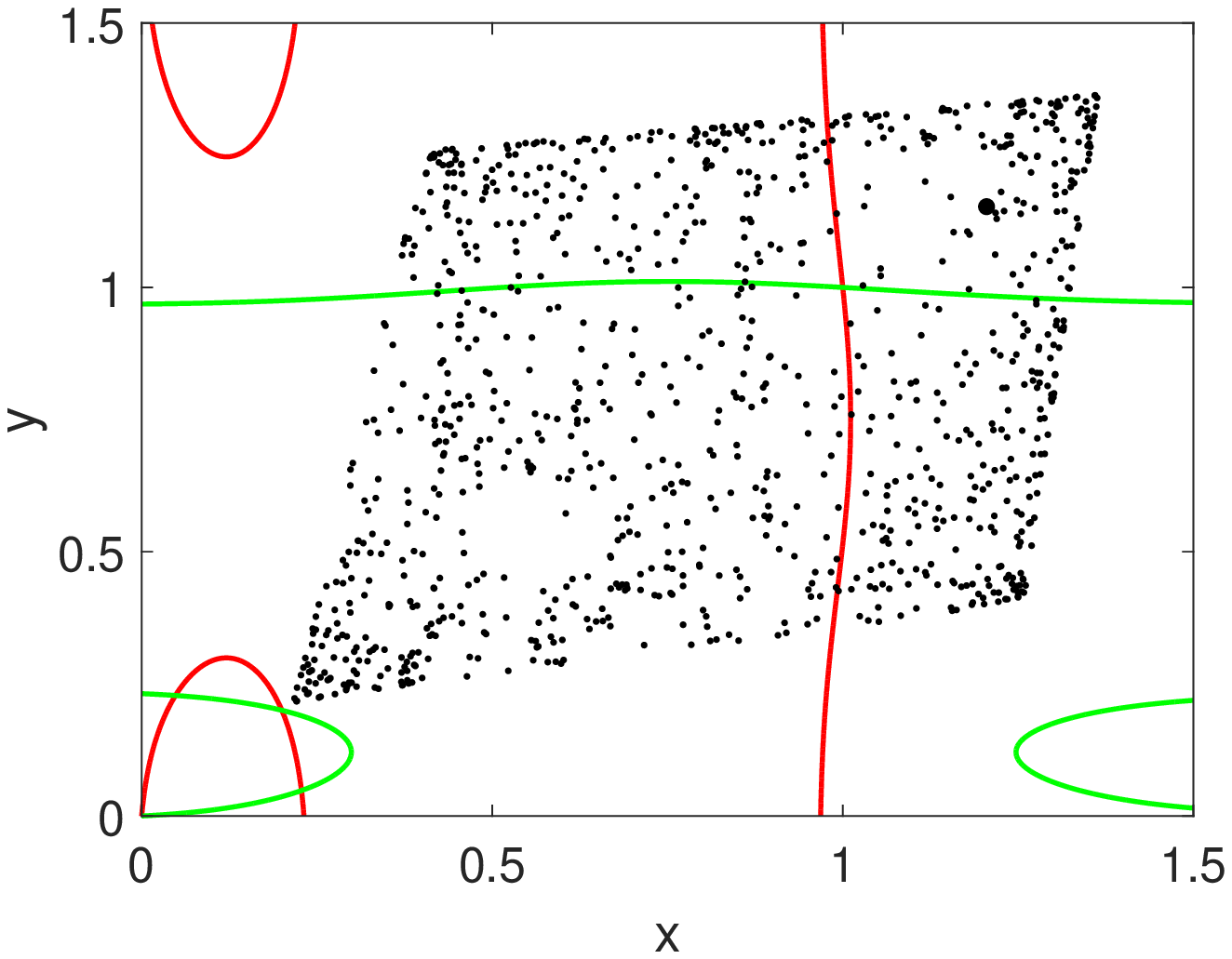}
\caption{$r=0.87$}
\end{subfigure}\hspace*{\fill}
\begin{subfigure}{0.5\textwidth}
\includegraphics[width=\linewidth]{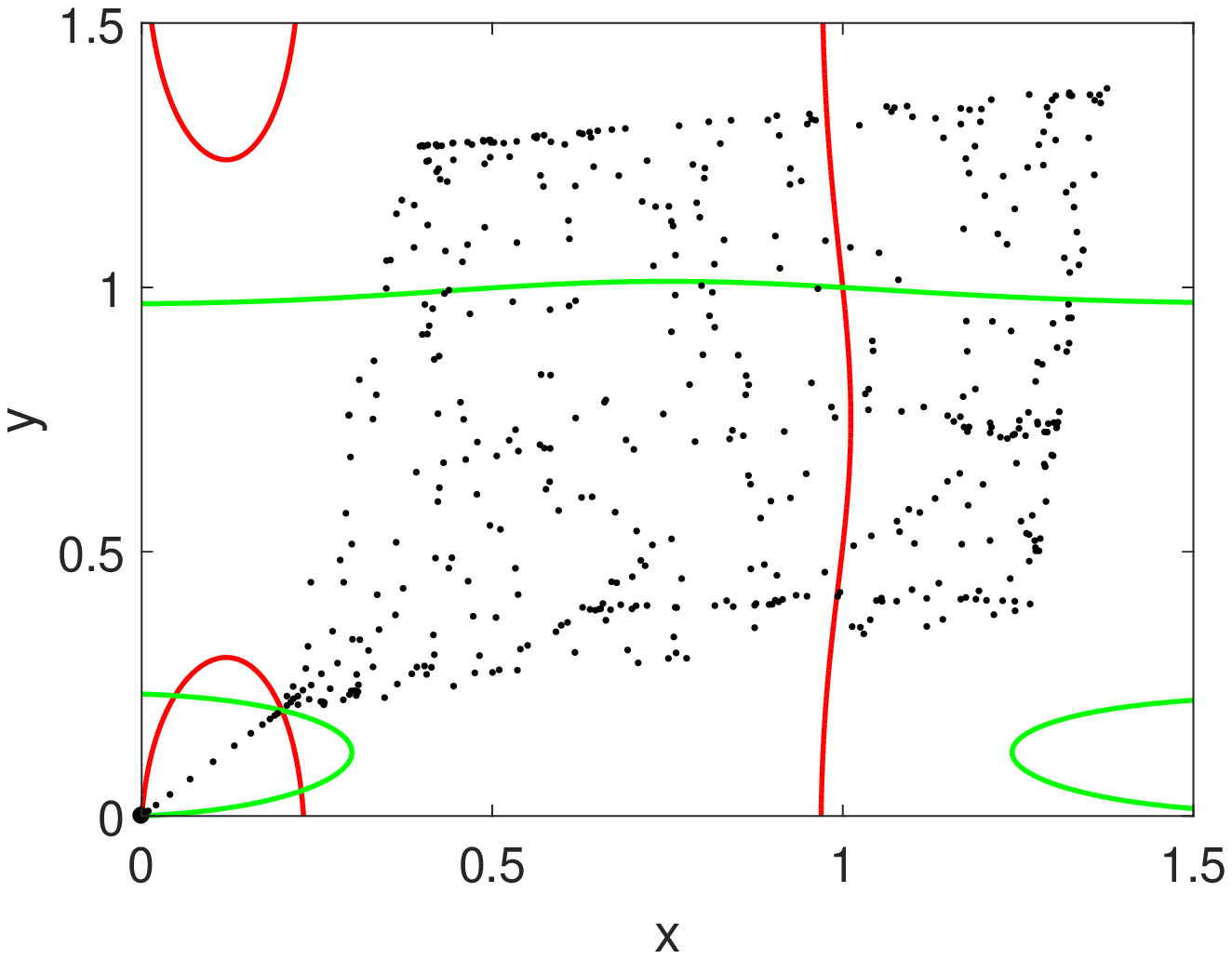}
\caption{$r=0.884$}
\end{subfigure}
\caption{Phase planes of the coupled system with $d=0.1$ and (a) $r=0.87$ and (b) $r=0.884$, between which a boundary crisis eliminates the symmetric coexistence attractor. Allee threshold $A = 0.2$, carrying capacity $K = 1$, $(A,A) < (x_0,y_0) < (\bar{A},\bar{A})$, 1000 time steps in both simulations, large symbols: final state. Nullclines in red and green, respectively.} \label{fig:test3}
\end{figure}

For $0.05<d<0.19$ and $d>0.28$, dispersal can not prevent essential extinction (see Figure \ref{fig:app}, $r>r_{th}$, outside of regions II and III, pink parameter region). The symmetric attractor is a chaotic rhombus (see Figure \ref{fig:test3}a) which disappears after a boundary crisis for $r>r_{th}$ and thus leads to essential extinction for almost all initial conditions (see Figure \ref{fig:test3}b). 

\subsection{Transients and crises} \label{Transients}
\begin{figure}
\begin{subfigure}{0.5\textwidth}
\includegraphics[width=\linewidth]{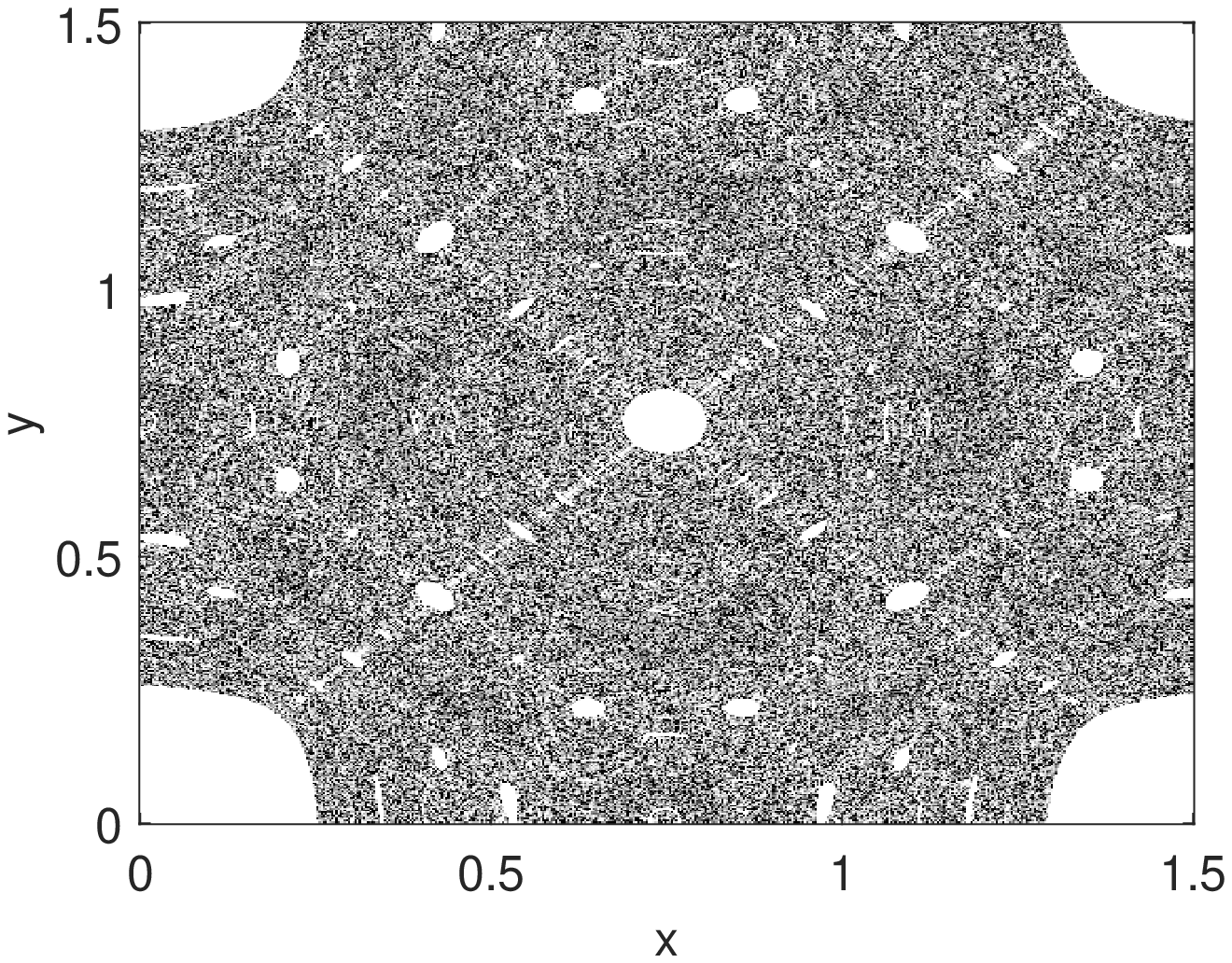}
\end{subfigure}\hspace*{\fill}
\begin{subfigure}{0.5\textwidth}
\includegraphics[width=\linewidth]{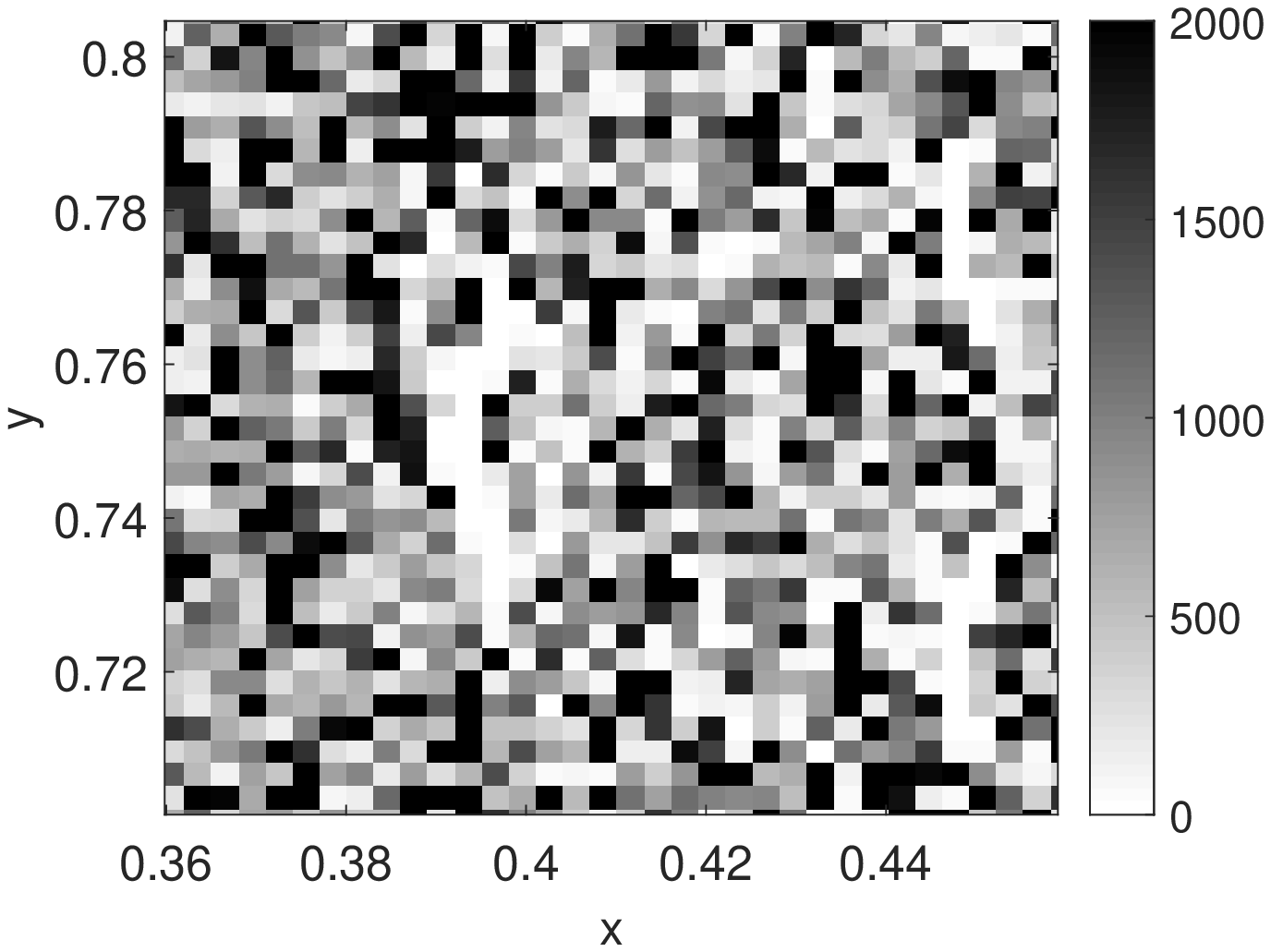}
\end{subfigure}
\caption{Left: time to extinction for parameter values $r=0.89$, $d=0.186$ and initial conditions $x_0$, $y_0 \in (0,1.5)$. Grey scale is chosen such that white means extinction after few time steps $t \approx 0$ and black means extinction at $t \approx 2000$ or later (see color bar). The population is called extinct at time $t$ when $x_t + y_t < 10^{-4}$. Right: enlarged section for selected initial conditions $x_0$, $y_0$.} \label{fig:time_to_ext}
\end{figure}
\begin{figure}
\centering
\includegraphics[width=0.7\linewidth]{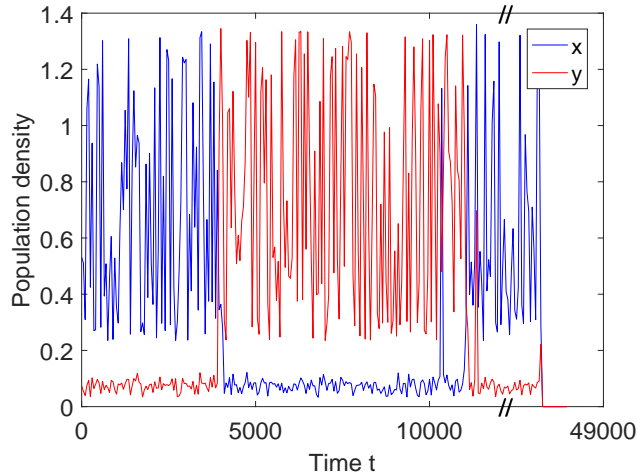}
\caption{Time series for parameter values $r=0.898$, $d=0.0415$ and initial conditions $x_0$, $y_0 = (0.07381,0.53102)$. Time steps $t \in (12000,45000)$ are hidden by a broken x-axis. Only every fifth value is plotted for better clarity.} \label{fig:long_transient}
\end{figure}

Transients are the part of the orbit from initial condition to the attractor and of particular importance in the case of crises \citep{hastings2018transient}. A boundary crisis occurs when an attractor exceeds the basin boundary around an invariant set, e.g. an equilibrium or a cycle \citep{neubert1997simple, vandermeer1999basin, wysham2008sudden, bischi2016qualitative,hastings2018transient}. Then the previous attractor forms a chaotic repeller or saddle and leads to long transients \citep{schreiber2003allee, wysham2008sudden}. \cite{schreiber2003allee} found long transients in a corresponding local model in parameter regions of essential extinction and proved that the time to extinction is sensitive to initial conditions due to the chaotic repeller formed by the basin boundary collision.\\
The transient behaviour which is seen in Figures \ref{fig:test2}, \ref{fig:test5} and \ref{fig:test3} can be partially explained with knowledge of the local system. We can also identify long transients induced by chaotic repellers or saddles. However, the coexistence of different persistence attractors can lead to different transient stages or transients that last orders of magnitudes longer than in the local case. In the following, we give numerical examples for both.\\
Different stages of transients before extinction of the population are shown in Figure \ref{fig:test5}. The approximately out-of-phase attractor (I) in Figure \ref{fig:test5}a undergoes a boundary crisis when $d$ decreases and merges with the transient chaotic rhombus that is also seen in Figure \ref{fig:test3}b. The two attractors disappear but are visible as ghosts (Figure \ref{fig:test5}b, I and II). Finally the population goes extinct (Figure \ref{fig:test5}b, III). In contrast to Figures \ref{fig:test2} and \ref{fig:test3}, the nullclines of the second iteration\footnote{The nullclines of the second iteration are the set of points satisfying $x_{t+2} = x_t$ for population $x$ and $y_{t+2} = y_t$ for population $y$, respectively \citep[cf.][]{kaplan2012understanding}.} in Figure \ref{fig:test5} highlight the invariant set at which the boundary crisis occurs (intersections of green and red nullclines). Figure \ref{fig:time_to_ext} presents the time to extinction for a range of initial population densities and the same parameters as used for Figure \ref{fig:test5}b. The sensitivity to initial conditions of transients is similar to the local system. The range of times until the population goes extinct reaches from values $\approx 0$ to more than $2000$ time steps (Figure \ref{fig:time_to_ext}). A steady-state analysis would not provide this information. From an ecological perspective, it is often more important to understand the transient than the asymptotic behaviour since this is on the relevant time scale. In contrast to regime shifts, where small parameter changes can lead to huge changes in the systems state, transient shifts can occur without additional environmental perturbations.\\
Figure \ref{fig:long_transient} shows a case of extremely long transients \citep{hastings2018transient}. The system passes the first 4700 time steps on one asymmetric ghost attractor until it switches to the other asymmetric ghost for the following 6000 time steps. Then the system switches back to the former ghost attractor, a behaviour that occurs due to a crisis in this parameter region. The long transient of about 34000 time steps ends abruptly and the population goes extinct after more than 46000 time steps without any parameter changes.
\section{Discussion and Conclusions} \label{Conclusion}
In this paper, we have developed a model for a spatially structured population with a local Allee effect and overcompensation. We found attractors to appear and disappear in the presence of dispersal. In contrast to \cite{knipl2016spatially} who state that the situation simplifies when dispersal increases, this conclusion does not hold for the model presented here. Nevertheless our results confirm two lines of research. Following \cite{amarasekare1998allee}, we showed that populations in patchy environments can have a large number of equilibria if both positive and negative density dependence are considered. We categorized extinction, symmetric and asymmetric attractors. Secondly, we identified additional symmetric attractors, analogous to \cite{hastings1993complex}. However, by Theorem \ref{th2} we gave conditions under which the behaviour of the coupled system can be derived from the behaviour of the uncoupled map. Overall, this simple model shows the complexity of interaction between chaotic dynamics, the Allee effect and dispersal.\\
In contrast to continuous-time models that suggest populations that are linked by dispersal to be more abundant and hence less susceptible to extinction \citep{amarasekare1998allee}, in discrete-time models not only small populations are endangered. However, we found two mechanisms that can prevent essential extinction of a spatially structured population whereas it takes place in the corresponding uncoupled system. Weak coupling of the two maps allows spatial asymmetry. Hence, it is possible to find one subpopulation with density above and one below the Allee threshold also for per-capita growth that leads to (essential) extinction without dispersal. Stronger coupling allows both subpopulations to persist above the Allee threshold due to (approximately) out-of-phase dynamics. Outside these parameter regions dispersal provides no mechanism to prevent essential extinction and the population goes extinct in almost all cases.\\ 
In summary, we support the conclusion of \cite{amarasekare1998allee} that interactions between Allee dynamics and dispersal create between-patch effects that lead to qualitative changes in the system. Populations are able to persist below the Allee threshold (rescue effect). Moreover, DIPEE provides another rescue effect for populations that suffer from essential extinction. The population with density below the Allee threshold is rescued from extinction and the population with density above the Allee threshold is rescued from essential extinction. Both subpopulations are prone to extinction without dispersal. However, a possibility for DIPEE is given only for specific initial conditions with a fractal basin boundary. For instance, DIPEE due to approximately out-of-phase dynamics for high dispersal benefits from asynchronous behaviour in the two patches \citep{lloyd1995coupled}. Small perturbations can synchronize this strongly connected system and thus lead to extinction \citep{earn2000coherence}.\\
Finally, we demonstrated the importance of the time scale since boundary crises may lead to long transients. Transient behaviour occurred also in the corresponding local system \citep{schreiber2003allee}. Our results for the coupled system support the statement that chaotic transients can last hundreds of time steps before the extinction state is reached. The duration of transients is also found to be sensitive to initial conditions. However, with the spatial structure of the model in this study, different persistence attractors can coexist. These can lead to different transient stages or transients that last orders of magnitudes longer than in the local case. A steady-state analysis will give no information about how long it takes a population to go extinct and what happens until extinction. On the other hand, short time series will eventually conceal that a population is damned to extinction for given parameters. Thus, a comprehensive analysis is fundamental to understand the complex behaviour of the presented system. This statement is supported for instance by \cite{wysham2008sudden} or \cite{hastings2018transient} who point out that ecologically relevant time scales are typically not the asymptotic time scales. In a next step, the impact of stochastic processes in the model could be tested since they are of particular importance in systems with multistability. Furthermore, a discrete-state model could be studied to investigate how lattice effects which inhibit chaos will lead to different dynamical behaviour \citep{henson2001lattice}. A question that we also do not address in this paper, is the significance of the chosen number of patches \citep{allen1993chaos, knipl2016spatially}. One could argue that in the case of more patches some effects may get lost or more pronounced. Further studies are needed to investigate the phenomena described (DIPEE, multiple attractors) on a broader spatial scale. Finally, the properties of dispersal could be refined in terms of asymmetric dispersal or dispersal mortality \citep{amarasekare1998allee, wu2020dispersal}.\\
Our model formulation is generic and does not depend on the Ricker growth model or the chosen implementation of the Allee effect. It is more about effects that are produced by coupled patches of locally overcompensatory dynamics with an Allee effect \citep{schreiber2003allee}. We tested other models of the same type and got similar results (not presented here). That is in line with \cite{amarasekare1998allee} and \cite{hastings1993complex}, who mention the generality of their results.\\ 
In summary, this paper contains some interesting results from the ecological and mathematical point of view: One key message is that small changes of parameters, perturbations or environmental conditions can have drastic consequences for a population. Even without external perturbations seemingly safe and unremarkable dynamics (long transients) can abruptly lead to extinction \citep{hastings2018transient}. This is of particular importance for species that show chaotic population dynamics. In this case they can be at risk not only for small population densities.\\
The effect of dispersal and connectivity can be either positive or negative. On the one hand dispersal can mediate local population persistence (rescue effect) or reduce overshoots and thus prevent essential extinction (DIPEE). On the other hand, dispersal can reduce local population sizes under the Allee threshold (Figure \ref{fig:app}, pink sprinkles in $r<r_{th}$) or induce an overshoot and thus cause (essential) extinction. These negative effects were not investigated in this work but should not be neglected. 

From the mathematical point of view it is interesting to observe a simple model setup with such a complexity in terms of multiple attractors and surprising results, e.g. long transients, caused by ghost attractors after various boundary crises \citep{hastings2018transient}.

 \section*{Acknowledgements} We thank Anastasiia Panchuk for providing the CompDTIMe routine for Matlab and two anonymous referees for providing very helpful comments which improved the manuscript. Alan Hastings and Sebastian Schreiber acknowledge support from the National Science Foundation under Grants DMS-1817124 and 1716803, respectively.

\section*{Appendix A}\label{appendixA}
\begin{theorem}\label{th1}
Let $f:[0,\infty) \rightarrow [0,\infty)$ be a three times continuous differentiable function that fulfills the following conditions:
\begin{itemize}
\item[(i)] $f$ has a unique critical point $D$ \vspace{0.1cm}
\item[(ii)] There exists an interval $[a,b]$ with \vspace{0.1cm}
\begin{itemize}
    \item[(a)] $f(x)>0 \text{ } \forall x \in [a,b]$ \vspace{0.1cm}
    \item[(b)] There is an $A \in (a,b)$ such that $f(A)=A$ and $f(x) \neq x \text{ } \forall x \in (0,A)$
    \item[(c)] $\lim\limits_{n \to \infty} f^n(x) =0 \text{ } \forall x \notin [a,b]$
    \item[(d)] The Schwartzian derivative of $f$ is negative for all $x \in [a,b]$
\end{itemize}
\end{itemize}
Define $A^*=max\{f^{-1}(A)\}$ and $M=f(D)$. Then: 
\begin{itemize}
\item \emph{Bistability:} If $f(M)>A$, then $f^n(x)\geq A \text{ } \forall n\geq0$, $x\in[A,A^*]$ and $\lim\limits_{n \to \infty} f^n(x) = 0 \text{ } \forall x\notin [A,A^*]$.\vspace{0.1cm}
\item \emph{Essential extinction:} If $f(M) <A$, then $\lim\limits_{n \to \infty} f^n(x)=0$ for Lebesgue almost every $x$.\vspace{0.1cm}
\item \emph{Chaotic semistability:} If $f(M) = A$, then the dynamics of $f$ restricted to $[A,A^*]$ are chaotic and $\lim\limits_{n \to \infty} f^n(x) = 0 \text{ } \forall x \notin [A,A^*]$.
\end{itemize}
\end{theorem}
According to \cite{schreiber2003allee}, we show that the criteria (i) and (ii) hold for function (\ref{f}) with parameter values $K=1$, $r>0$ and $0<A<K$. That is, (\ref{f}) shows either bistability or essential extinction, depending on the parameter values.\\
\begin{itemize}
 \item[(i)] $f$ has a unique positive critical point $D$ at:
\begin{align*}
x = \frac{1+A}{4}+\frac{1}{4} \sqrt{\frac{8A+r+2Ar+A^2r}{r}}
\end{align*}
\item[(ii)] For an interval $[a,b]$ where $a,b>0$ 
\begin{enumerate}
\item[(a)] is fulfilled by the product of two positive values ($x$ itself and the exponential function).
\item[(b)] is fulfilled by the Allee threshold $A$. For all $0 < x < A$ we get $f(x)<x$ since the exponential function has a negative exponent and is thus smaller than one.
\item[(c)] Choose $a \in (0,A)$. $\lim\limits_{x \to \infty} f(x) =0$ and $f$ has a unique positive critical point. Thus, there exists a unique $b>a$ such that $f(b) = a$. It follows that $f(x) \in [0,A]$ for all $x \in [b,\infty)$. Hence, (c) is fulfilled.
\item[(d)] The Schwartzian derivative of $f$ is:
\begin{align*}
Sf(x)=-\frac{q_1^2(x)q_2(x)+\frac{12r^2x^2}{A^2}+\frac{12r}{A} }{2(1+xq_1)^2}
\end{align*}
with
\begin{align*}
q_1(x) &= r\left(\frac{1-2x}{A}+1\right)\\
q_2(x) &= 6+x^2+q_1^2(x)+4xq_1(x)
\end{align*}
All terms except $q_2(x)$ are obviously positive. For $q_2(x)$, we have:
\begin{itemize}
\item $q_2(0) = 6$
\item The only minimum (for positive values of $x$) occurs at 
\begin{align*}
x_{min}=\frac{1+A}{4}+\frac{1}{4} \sqrt{\frac{16A+r+2Ar+A^2r}{r}}
\end{align*}
with $f(x_{min}) = 2$. 
\item $\lim\limits_{x \to \infty} q_2(x) = \infty$.
\end{itemize}
In summary, $q_2$ is also positive. Thus, the Schwartzian derivative is negative for all $r>0$, $A>0$ and $x>0$. 
\end{enumerate}
\end{itemize}

\section*{Appendix B: Proof of Theorem~\ref{th2}}
To prove the theorem, let $\widetilde F(x,y)=(f(x),f(y))$ denote the uncoupled map and let $F(x,y)=((1-d)f(x)+df(y),df(x)+(1-d)f(y))$ be the coupled map with $d>0$. Assume that $f$ has a linearly stable periodic orbit $\mathcal{O}=\{p,f(p),\dots,f^{n-1}(p)\}$ of period $n$ with $p\in[A,\infty)$. Since $f$ has a negative Schwartzian derivative and a single critical point on the interval $[A,\infty)$, Theorem A of \citet{vanstrien-81} implies that the complement of the basin of attraction of ${\mathcal O}$  for $f$ can be decomposed into a finite number of compact, $f$-invariant sets which have a dense orbit and are hyperbolic repellers: there exists $c>0$ and $\lambda>1$ such that $|(f^t)'(x)|\ge c \lambda^t$ for all points $x$ in the set and $t\ge 1$. Consequently, the $2$-dimensional mapping $\tilde F$ is an Axiom A endomorphism~\citep[pg. 271]{przytycki-76}: the derivative of $\widetilde F$ is non-singular for all points in the non-wandering set $\Omega(\widetilde F)=\{(x,y)\in \C:$ for every neighborhood $U$ of $(x,y)$, $\widetilde F^t(U) \cap U \neq \emptyset$ for some $t\ge1\}$, $\Omega(\widetilde F)$ is a hyperbolic set, and the periodic points are dense in $\Omega (\widetilde F)$. \citet[3.11-3.14 and 3.17]{przytycki-76} imply that (i) $\Omega(\widetilde F)$ decomposes in a finite number of invariant sets $\Omega^1(\widetilde F),\dots, \Omega^m(\widetilde F)$ and (ii) maps sufficiently $C^1$ close to $\widetilde F$ are Axiom A endomorphisms whose invariant sets $\Omega^i(F)$ are close to $\Omega^i(\widetilde F)$. In particular, property (ii) implies that $F(x)$ is an Axiom A endomorphism provided that $d>0$ is sufficiently small. The invariant sets $\Omega^i(\tilde F)$ for $\tilde F$ correspond the linearly stable periodic orbits defined by $\mathcal{P}$, and products of the hyperbolic repellers for $f$ and the linearly stable periodic orbits of $f$. Without loss of generality, let $\Omega^i(\tilde F)$ for $1\le i\le n+3$ correspond to the linearly stable periodic orbits of $\tilde F$ and $\Omega^i(\tilde F)$ for $i>n+3$ correspond to the saddles and repellers of $\tilde F$. For $d>0$ sufficiently small, $\Omega^i(F)$ retain these properties. For $d\ge 0$ sufficiently small, the proof of Theorem IV.1.2 in \citet{qian_xie2009} implies that the complement of the basin attraction of $\cup_{i=1}^{n+3}\Omega^i(F)$ has Lebesgue measure zero.\\
To prove assertion (iii), assume $n$ is not a power of $2$. Then \citet{sharkovsky1964coexistence} proved that $f$ has an infinite number of periodic orbits. All but two of these periodic orbits lie in the hyperbolic repellers of $f$. Consequently, the set of saddles and repellers $\cup_{i>n+3}\Omega^i(\tilde F)$ of $\tilde F$ contain an infinite number of periodic points. Hyperbolicity of these saddles and repellers implies that the set of saddles and repellers $\cup_{i>n+3}\Omega^i(F)$ of $F$ has an infinite number of periodic orbits for $d>0$ sufficiently small.

\medskip
\bibliographystyle{spbasic}      

\bibliography{mybib}

\end{document}